\documentclass[12pt]{article}
\usepackage{latexsym}
\usepackage{graphicx}
\usepackage{verbatim}
\usepackage{amssymb}
\usepackage{amsmath}
\title{Generalised dimensions of measures on almost self-affine sets}
\author{K.J. Falconer\\
\small{{\it Mathematical Institute,  
University of St~Andrews, North Haugh, St~Andrews,}} \\
\small{{\it Fife, KY16~9SS, Scotland }}} 
 
\def\bbbr{{\mathbb R}}
\def\bbbn{{\mathbb N}}

\newtheorem{theo}{Theorem}
\newtheorem{prop}[theo]{Proposition}

\newtheorem{lem}[theo]{Lemma}
\newtheorem{cor}[theo]{Corollary}

\newcommand\hdd{\mbox{\rm dim}_{\rm H}\,} 
\newcommand\ubd{\overline{\mbox{\rm dim}}_{\rm B}\,} 
\newcommand\lbd{\underline{\mbox{\rm dim}}_{\rm B}\,} 
\newcommand\qd{D_{q}\,}
\newcommand\lqd{\underline{D}_{q}\,}
\newcommand\uqd{\overline{D}_{q}\,}

\newcommand\Esp{\bbbr^{N}}

\newcommand{\svs}[1]{\phi^{s}({#1})} 
\newcommand{\svf}[2]{\phi^{{#1}}({#2})} 
\newcommand{\ssi}{_{\bf i}} 
\newcommand{\ti}{T_{\bf i}} 
\newcommand{\bi}{{\bf i}} 
\newcommand{\bj}{{\bf j}} 
\newcommand{\bu}{{\bf u}} 
\newcommand{\bv}{{\bf v}} 
\newcommand{\bw}{{\bf w}} 
\newcommand{\ii}{I_{\infty}} 
\newcommand{\ci}{C_{\bf i}} 
\newcommand{\w}{\omega} 
\newcommand{\xw}{\Pi^{\omega}} 
\newcommand{\E}{{\sf E}} 
\renewcommand{\P}{{\sf P}} 
\newcommand{\J}{{\cal J}} 
\newcommand{\ps}{\phi^{s}} 
\newcommand{\be}{\begin{equation}} 
\newcommand{\ee}{\end{equation}} 

\begin{document}
\maketitle

\begin{abstract} 
We establish a generic formula for the generalised $q$-dimensions of measures 
supported by almost self-affine sets, for all $q>1$.  These $q$-dimensions 
may exhibit phase transitions as $q$ varies. We first consider general measures and then specialise to Bernoulli and Gibbs measures. Our method involves estimating expectations of moment expressions in terms of `multienergy' integrals which we then bound using induction on families of trees.  
\medskip

\noindent AMS classification scheme numbers: 28A80, 37C45\end{abstract}

\section{Introduction}
\setcounter{equation}{0}

Let  $S_{1},\ldots ,S_{m}:\bbbr^{N}\rightarrow \bbbr^{N}$ 
be a family of  contractions; 
thus there are constants $c_{i}<1$ such that $|S_{i}(x) - S_{i}(y)| 
\leq c_{i}|x-y|$  for all 
$x,y \in \bbbr^{N}$.  Such a family is 
known as an {\it iterated function system} (IFS), and it is 
well-known that there exists a unique non-empty compact set $E$ 
satisfying
$$ E= \bigcup_{i=1}^{m} S_{i}(E),$$ 
called the {\it attractor} of the system.  
If the $S_{i} = T_{i} + a_{i} \,(1=1,\ldots,m)$ are affine contractions, where $T_{1},\ldots ,T_{m}$ are
non-singular contracting linear mappings on $\Esp$ and $a_{1},\ldots 
,a_{m} \in \Esp$ are translation vectors, we call $E$  a {\it self-affine set}, see \cite{Fa}. 

The Hausdorff and box-counting dimensions of many self-affine sets $E$ are given by  $\min\{d(T_{1},\ldots, T_{m}) , N \}$, where 
$$d(T_{1},\ldots, T_{m}) = \inf \Big\{s: \sum_{k=1}^{\infty}
\sum_{\bi \in I_{k}}\svs {T_{\bi}} < \infty\Big \},$$
a number sometimes called the {\em affinity dimension} of $E$, which is defined in terms of the singular values  of iterated compositions of the mappings $T_i$, see Section 2 for details. 
This was shown in \cite{Fa1} to hold for almost all $(a_{1},\ldots ,a_{m}) \in \bbbr^{mN}$
provided that $\| T_{i} \| < \frac{1}{3}$ for all $i$, a restriction that was soon weakened to $\| T_{i} \| < \frac{1}{2}$ in  \cite{Sol}.
The affinity dimension turns out to give the Hausdorff or box dimensions of self-affine sets in many other cases, see for example \cite{Fa2,HL,KS}, and is regarded as a `generic' formula for the dimensions of self-affine sets. Nevertheless,  the dimensions of self-affine sets need not vary continuously with 
the parameters $(a_{1},\ldots ,a_{m})$ and for highly regular 
constructions such as self-affine carpets (where the $S_{i}$ map a 
square onto rectangles selected from a rectangular grid) the dimension 
of the self-affine set is in general strictly less than than its affinity dimension, see \cite{Bed, McM}. For more on the dimensions of self-affine sets, see the surveys \cite{Edg, PSol}. 

Recently, Jordan,  Pollicott and Simon \cite{JPS} introduced a variant of self-affine sets, which we will term {\it almost self-affine} sets, with rather more randomness allowed in the translation parameters than for strictly self-affine sets. Here an independent random perturbation is made at each stage of the iterated construction of the set, yielding a `statistically self-affine' set, and this was shown to have Hausdorff and box  dimensions $\min\{d(T_{1},\ldots, T_{m}) , N \}$ almost surely, with the only restriction that $\| T_{i} \| < 1$ for all $i$; see (\ref{2.6}) below. (Note that it is often convenient to use the language of probability rather than of measure theory when considering such constructions.)

It is natural to consider multifractal analogues of these dimension formulae, in particular to seek the
 generalised $q$-dimensions (also termed generalised R\'{e}nyi 
dimensions) of measures supported by 
self-affine sets.  The generalised $q$-dimension of a measure $\tau$ reflects the behaviour of moment sums of $\tau$ over small boxes; it 
is given for $q> 0, q \neq 1$ by  
$$D_{q}(\tau) = \lim_{r \rightarrow 0}\frac{\log  M_{r}(q)}
{(q-1) \log r}$$
provided the limit exists (or taking lower and upper values otherwise), where
$M_{r}(q) = \sum_{C \in {\cal M}_{r}}\tau (C)^{q}$, with the sum 
over the set ${\cal M}_{r}$ of  
mesh cubes of side $r$, see Section 3 for more details.  

There is a natural analogue of the affinity dimension that is appropriate for  $q$-dimensions of measures supported by self-affine sets, namely (for $q>1$)
\be
d_{q}^-(T_{1},\ldots, T_{m};\mu) =
\sup \{s :\sum_{k=0}^{\infty} \sum_{\bi \in I_{k}}  
 \svs {T_{\bi}} ^{1-q} 
\mu (C_{\bi})^{q} <\infty \},\label{affdim}
\ee
where $\mu$ is a measure on the underlying code space (for $0<q<1$ the supremum is replaced by an infimum). One might hope that the (lower) generalised dimensions of self-affine measures would equal $\min\{d_q^-, N\}$  in a `generic' sense. Again, examples such as measures on self-affine carpets show that this cannot be true for all constructions, see \cite{BM,K,O}.
Nevertheless, we showed in \cite{Fa5} that  $\min\{d_q^-, N\}$ gives an upper bound  if $q>0$, and in the case of $1<q \leq 2$ equals the generalised dimension for almost all  sets of translation vectors $(a_{1},\ldots ,a_{m}) \in \bbbr^{mN}$, provided that  $\| T_{i} \| < \frac{1}{2}$ for all $i$.

To estimate higher moments, that is to find $D_{q}$ for $q >2$, more randomness seems to be needed, and a very natural setting is for measures supported by almost self-affine sets. Our main result, Theorem \ref{thm8.1}, is that the (lower) generalised $q$-dimension of a measure on an almost self-affine set  equals 
$\min\{d_q^-, N\}$ almost surely for all $q>1$.  We first consider general measures on almost self-affine sets, and then specialise to Bernoulli measures and Gibbs measures in Corollaries \ref{cor8.3} and \ref{cor8.5} As with self-affine measures for $1<q \leq 2$, see \cite{Fa5},  the generalised $q$-dimensions can exhibit phase transitions, corresponding to the non-differentiability of $d_q^-$ at those $q$ where where $d_q^-$ is an integer.

Upper bounds for the  $q$-dimensions follow from routine methods. Obtaining almost sure lower bounds for the  $q$-dimensions is much more involved and breaks into two stages. First, in Section 6, we show that the expectation of (an equivalent integral version of)  $M_{r}(q)$ is controlled by certain `multienegy integrals'  of the form$$\int \cdots \int 
 \phi^s(\bi_1,\ldots,\bi_q)^{-1} d\mu(\bi_1)\ldots d\mu(\bi_q)$$
  (in the special case where q is an integer)
where $\phi^s(\bi_1,\ldots,\bi_q)$ is given in terms of products of singular value functions of iterated products of the $T_i$.  Then in Section 7 we bound these integrals by breaking up the domain of integration and estimating the integral over each such subdomain using induction on families of trees. This leads to the desired almost sure lower bounds for the $q$-dimensions in Section 8.
 
\section{Definitions and notation}
\setcounter{equation}{0}
\setcounter{theo}{0}

We will work throughout with  contracting,
non-singular, linear mappings $T \in {\cal L}(\bbbr^{N},\bbbr^{N})$; of course  products of such 
mappings will also be contracting and non-singular.  Recall that the {\it singular values} 
$\alpha_{i} \equiv  \alpha_{i}(T)$ of $T \; (i = 1,\ldots ,N)$ are the positive square roots 
of the eigenvalues of $TT^{\ast}$, 
where $T^{\ast}$ is the transpose of $T$, or  
equivalently are the lengths of the (mutually 
perpendicular) principal semiaxes of $T(B)$, where $B$ is the unit ball 
in $\bbbr^{N}$.  We adopt the convention that
$1> \alpha_{1} \geq \alpha_{2} \geq \ldots \geq \alpha_{N}>0$.
The {\it singular value function } 
$\phi^{s}(T)$  is 
central in the analysis of self-affine sets.  For $0 \leq s \leq N$ we define
$$
\phi^{s}(T) = \alpha_{1} 
\alpha_{2} \ldots \alpha_{j-1}\alpha_{j}^{s-j+1}, 
$$
where $j$ is the 
integer such that $j-1 < s \leq j$. It is convenient to 
set $\phi^{s}(T)= (\alpha_{1} \alpha_{2} \ldots \alpha_{N})^{s/N} = 
(\det T)^{s/N}$ for $s>N$.  

Clearly $\svs T$ is continuous and 
strictly decreasing in $s$, and
 is  
{\it sub-multiplicative}, that is, for all $s \geq 0$,
\be
\svs {TU} \leq  \svs T \svs U  \label{2.2}
\ee
for all $T,U \in {\cal L}(\bbbr^{N},\bbbr^{N})$, see \cite{Fa1} for these basic properties.

Many fractals, including self-similar, self-affine and almost self-affine sets may be constructed in a hierarchical manner 
which can conveniently be indexed by a {\it code space} or {\it sequence space}.
For $k=0,1,2, \ldots$ let $I_{k}$ be the
set of all $k$-term sequences or words formed from the integers $1,2, \ldots, m$,
that is $I_{k} = \{ (i_{1}, i_{2}, \ldots , i_{k}): \, 1 \leq i_{j} \leq 
m \}$;
we take $I_{0}$ to contain just the empty word $\emptyset$. 
 We often abbreviate a
word in $I_{k}$ by
${\bf i} = (i_{1}, i_{2}, \ldots , i_{k} ) $ and write $|{\bf i}|=k$ for 
the length of ${\bf i}$.
We write
$I= \cup^{\infty}_{k=0} I_{k} $
for the set of all such finite words, and $I_{\infty}$ for the corresponding
set of infinite words, so
$I_{\infty} = \{(i_{1}, i_{2}, \ldots ): 1 \leq i_{j} \leq m \}$.
Juxtaposition
of ${\bf i}$ and ${\bf j}$ is written ${\bf ij}$.  
We write ${\bf 
i}|_{k} = (i_{1}, \ldots , i_{k})$ for the {\it curtailment} after $k$ terms
of ${\bf i} = (i_{1}, i_{2}, \ldots ) \in I_{\infty}$, 
or of   ${\bf i} = (i_{1}, \ldots,  i_{k'}) \in I$ if $ k 
\leq k^{\prime}$.  We write ${\bf i} \preceq {\bf j}$ to mean that ${\bf i}$ is a 
curtailment of ${\bf j}$.  If ${\bf i,j} \in I_{\infty} $  then 
${\bf i}\wedge {\bf j}$ is the maximal sequence such that both 
${\bf i}\wedge {\bf j} \preceq {\bf i}$ and ${\bf i}\wedge {\bf j} \preceq {\bf j}$.

We may topologise $I_{\infty}$ in a natural way by the metric 
$d(\bi,\bj) = 2^{-|\bi \wedge \bj |}$ for distinct $\bi,\bj \in 
I_{\infty}$ which makes $I_{\infty}$  into a compact metric space. 
The {\it cylinders}  
$C \ssi = \{\bj \in I_{\infty} : \bi \preceq \bj \}$ for  $\bi \in I$ form a base 
of open and closed neighbourhoods of $I_{\infty}$.

 It is convenient to identify $I$ with the vertices of an $m$-ary
 rooted tree with root $\emptyset$. The edges of this tree join each vertex $\bi \in I$  to its $m$  `children'  
 $\bi_1,\ldots,\bi_m$. The estimates in Section \ref{integralest} 
 involve certain automorphisms of this tree.

Compositions of the contractions $T_{1},\ldots ,T_{m}$ will be written
$T\ssi \equiv T_{i_{1}}T_{i_{2}} \ldots T_{i_{k}}$ where ${\bf i} = (i_{1}, i_{2}, \ldots , 
i_{k} ) $, with 
$T_{\emptyset}$ the identity mapping.  
Set 
\begin{eqnarray}
a_{-} & = & \min_{1 \leq i \leq m} \alpha_{n}(T_{i}) 
\label{2.4c1} \\
a_{+} & = & \max_{1 \leq i \leq m} \alpha_{1}(T_{i}) \label{2.4c}
\end{eqnarray}
where $\alpha_{j}(T_{i})$ 
are the singular values of $T_{i}$.  Then $0 < a_{-} \leq a_{+} <1 $, 
and 
\begin{equation*}
a_{-}^{|\bi|} \leq \alpha_{j}(T_{\bi}) \leq a_{+}^{|\bi|}
\end{equation*}
 for all $\bi \in I$ and $j= 1,\ldots,n$, so that
\begin{equation}
a_{-}^{s|\bi|} \leq \svs {T_{\bi}} \leq a_{+}^{s|\bi|}.
\end{equation}
We also note that, for $h>0$, 
\begin{equation}
a_{-}^{h|\bi|}\svs {T_{\bi}}  \leq \phi^{s+h}(T_{\bi})\leq a_{+}^{h|\bi|}\svs {T_{\bi}}. \label{phicont}
\end{equation}

We now introduce a notation that will permit a random perturbation at each stage of the 
hierarchical construction of the attractor as in \cite{JPS}.
Let $D$ be a bounded region of $\Esp$ and for each $\bi \in I$ 
let $\w_{\bi} \equiv \w_{i_1,\ldots,i_k} \in D$ be a `displacement' or   `perturbation'  which will eventually be random. Let $\w = \{\w_\bi: \bi \in I\}$ denote the aggregate of the $\w_{\bi}$.
Define the projection $\xw: I_{\infty} \to \Esp$ by
\begin{eqnarray}
\xw(\bi) & =&\lim_{k \to \infty} (T_{i_1} + \w_{i_1})  (T_{i_2} + \w_{i_1,i_2}) \cdots
  (T_{i_k} + \w_{i_1,\ldots,i_k})(x)\label{points1}\\
& =&  \w_{i_1} +T_{i_1}\w_{i_1,i_2} + T_{i_1}T_{i_2}\w_{i_1,i_2,i_3} 
+ \cdots. \label{points2}
\end{eqnarray}
It is easily checked that this limit exists and is independent of $x \in \Esp$, and that the map
$\bi \mapsto \xw(\bi)$ is continuous.

We term the compact set 
\begin{equation}
E^{ \w} = \bigcup_{{\bf i} \in I_{\infty}} \xw(\bi) \subseteq \Esp,   \label{2.6}
\end{equation}
 an {\it almost self-affine set}.  Note that if $B \subseteq \Esp$ is a ball large enough so that  $T_{i_1} (B)+ \w_{i_1,\dots,i_k} \subseteq B$ for all $i_1,\dots,i_k$ then
\be
E^{ \w} = \bigcap_{k=0}^{\infty} 
\bigcup_{i_1,\dots,i_k \in I_k}  (T_{i_1} + \w_{i_1}) (T_{i_2} + \w_{i_1,i_2}) \cdots (T_{i_k} + \w_{i_1,\ldots,i_k})(B), \label{construct}
\ee
which represents the standard hierachical way of constructing $E^{ \w}$.

A standard covering argument, involving 
dividing up each of the sets in (\ref{construct}) into appropriate pieces, shows that for all $\w$, the Hausdorff and lower and upper box-counting dimensions satisfy
\be
\hdd (E^{ \w} ) \leq \lbd (E^{ \w} ) \leq \ubd (E^{ \w} ) \leq \min\{d(T_{1},\ldots, T_{m}),N\} ,\label{setdim}
\ee
where
\be
d(T_{1},\ldots, T_{m}) = \inf \Big\{s: \sum_{k=1}^{\infty}
\sum_{\bi \in I_{k}}\svs {T_{\bi}} < \infty\Big \};\label{dimexp}
\ee  
see \cite{Fa1,JPS}. 

There are many situations where equality holds in (\ref{setdim}). These include self-affine sets  (where  $\w_{i_1,\ldots,i_k} = \w_{i_1}$ depends only on the first subscript) for ${\cal L}^{mN}$-almost all $\w_{i_1},\ldots, \w_{i_m}$, provided
$\|T_i\| < \frac{1}{2}$ for all $i$, see \cite{Fa1, Sol}. Equality also holds with probability one for the almost self-affine sets introduced in \cite{JPS}, where the $\{\w_\bi: \bi \in I\}$  are independent with identical distributions of bounded density within a region $D$. 

We now introduce measures supported on $E^{ \w}$ by projecting a measure from $I_\infty$.
Let $\mu$ be a finite Borel measure (with respect to the 
metric $d$) on $\ii$.  For each $\w = \{\w_\bi: \bi \in I\}$ let $\mu^{\omega}$
 be the image of $\mu$ under the projection $\xw$, that is
\begin{equation}
\mu^{\omega}(A) = \mu \{ \bi : \xw(\bi) \in A \} \label{mesdef}
\end{equation}
for $A \subseteq \bbbr^{N}$, or equivalently by 
\begin{equation}
\int f (x) d\mu^{\omega}(x) = \int f (\xw(\bi)) d\mu(\bi) \label{mesintdef}
\end{equation}
for continuous $f: \bbbr^{N} \rightarrow \bbbr$. For each $\w$ the measure $\mu^{\omega}$ is supported by the almost self-affine set  $E^{ \w}$. If $\w_{i_1,\ldots,i_k} = w_{i_1}$ for all $\bi = i_1,\ldots,i_k$ we get the self-affine measures studied in \cite{Fa5}.

\section{Generalised $q$-dimensions}
\setcounter{equation}{0}
\setcounter{theo}{0}

One approach to multifractal analysis of 
a measure on $\Esp$ involves generalised 
$q$-di\-men\-sions; see \cite{Fa3,Gr,Har,Man,P} 
for various treatments.  
The generalised dimensions of a finite Borel measure $\tau$ of bounded 
support may be defined along the lines of
box-counting dimension using {\it $r$-mesh cubes}, that is cubes in 
$\Esp$ of the form $[j_{1}r,(j_{1}+1)r) \times \cdots \times
[j_{n}r,(j_{n}+1)r)$ where $j_{1}, \ldots , j_{n}$ are integers.  
We write  ${\cal M}_{r}$ for the set of $r$-mesh cubes in 
$\bbbr^{N}$.  
The $q$-dimensions reflect the power law 
behaviour of moment sums of $\tau$. 
For $q >0$ and $r>0$ set
\begin{equation}
M_{r}(q) = \sum_{C\in  {\cal M}_{r}} \tau (C)^{q},\label{mom}
\end{equation}
where the sum is over the $r$-mesh cubes $C$ such that $\tau (C)>0.$ 
We identify the power law behaviour of $M_{r}(q)$ 
by defining, for  $q\neq 1$, the {\it lower} and {\it upper generalised q-dimensions} 
of $\tau$ 
\be
\lqd (\tau) = \liminf_{r \rightarrow 0}\frac{\log  M_{r}(q)}{(q-1) \log r}
 \quad \mbox{and} \quad 
\uqd (\tau) = \limsup_{r \rightarrow 0}\frac{\log  M_{r}(q)}{(q-1) \log 
r}.\label{3.b} 
\ee
If, as frequently happens, $\lqd (\tau)=\uqd (\tau)$, we write $D_{q}(\tau)$ for the common 
value which we refer to as the {\it generalised q-dimension}.

It is easily verified that $\lqd (\tau)$ and $\uqd (\tau)$ are each 
nonincreasing in $q$ and continuous (for  $q\neq 1$), and that 
$0 \leq \lqd (\tau) \leq \uqd (\tau) \leq N$ for all $q$. 

In this paper we will be entirely concerned with  higher moments and  will assume that $q>1$ throughout. In this case, the definitions of $q$ dimensions are independent of the origin 
and orientation chosen for the mesh cubes. 

There are useful integral forms of  
$\lqd$ and $\uqd$.  For $q > 1$, 
\begin{eqnarray}
\lqd (\tau) &=& \liminf_{r \rightarrow 0}\frac{\log \int 
\tau(B(x,r))^{q-1}d\tau(x)}{(q-1) \log r} \label{intdefl}\\
\mbox{and} \quad 
\uqd (\tau) &=& \limsup_{r \rightarrow 0}\frac{\log \int 
\tau(B(x,r))^{q-1}d\tau(x)}{(q-1) \log r},\label{intdef} 
\end{eqnarray}
see \cite{Lau}. 

\section{Upper bounds for generalised dimensions}
\setcounter{equation}{0}
\setcounter{theo}{0}

It is not difficult to derive natural upper bounds for $\lqd (\mu^\omega)$ and $\uqd 
(\mu^\omega)$ valid for a general measure $\mu$ and 
all $\omega$.  For given $s> 0$
and $0<r<1$ let $j$ be the integer 
such that $j-1 < s \leq j$ and define 
\begin{eqnarray}
J^{s}(r)= \{\bi =(i_{1},\ldots,i_{k}) \in I: 
\alpha_{j}(T_{i_{1},\ldots,i_{k}}) \leq r 
<\alpha_{j}(T_{i_{1},\ldots,i_{k-1}}) \}. \label{Jr}
\end{eqnarray} 
The finite set of sequences $J^{s}(r)$ is a {\it cut-set} 
or {\it stopping} in the sense that for every $\bi \in I_{\infty}$ there 
is a unique integer $k$ such that $\bi |_k \in J^{s}(r)$. 
From (\ref{2.4c1})
$
a_{-}r < \alpha_{j}(T_{\bi}) \leq r$ for all $\bi 
\in J_{r}.$ The basic estimate is as follows.

\begin{prop}\label{sumest}
Let $\mu$ be a finite Borel measure on $\ii$, let $\mu^{\omega}$ be the measure on $E^\w$ 
defined by {\rm (\ref{mesdef})}. For $q \geq 1$ and $0<s\leq N$
there is a number $c>0$ such that, for 
all $\w$ and all sufficiently small $r$,
\begin{eqnarray}
r^{s(1-q)} \sum_{C \in {\cal M}_{r}} \mu^{\w}(C)^{q}
\geq c \sum_{\bi \in J^{s}(r)}
\svs {T_{\bi}}^{1-q}\mu(C_{\bi})^{q}. \label{3.2n}
\end{eqnarray} 
\end{prop} 

\noindent{\it Proof.}  
The proof in \cite[Proposition 4.1]{Fa1} holds virtually unchanged, by covering the ellipsoids  
$ (T_{i_1} + \w_{i_1}) (T_{i_2} + \w_{i_1,i_2}) \cdots  (T_{i_k} + \w_{i_1,\ldots,i_k})(B)$ for $(i_{1},\ldots,i_{k}) \in J^s(r)$
by cubes of sidelengths $\alpha_j(T_{i_1}T_{i_2}\cdots T_{i_k})$, where $j$ is the integer 
such that $j-1 < s \leq j$, and summing the measures of these cubes and using Jensen's and Minkowski's inequalities. 
$\Box$

\medskip

This leads us to define quantities that one might hope would give the lower and upper generalised $q$-dimensions for $q>1$.
\begin{eqnarray} 
d_{q}^{-} &\equiv&  d_{q}^{-}(T_{1},\ldots,T_{m};\mu) = 
\sup \{s : \limsup_{r\rightarrow 0}\sum_{\bi \in J^{s}(r)} \svs 
{T_{\bi}}^{1-q}\mu(C_{\bi})^{q} < \infty \}, \label{basicb1}\\
d_{q}^{+} &\equiv&  d_{q}^{+}(T_{1},\ldots,T_{m};\mu) = 
\sup \{s : \liminf_{r\rightarrow 0}\sum_{\bi \in J^{s}(r)} \svs  
{T_{\bi}}^{1-q}\mu(C_{\bi})^{q} < \infty \}, \label{basicb2}
\end{eqnarray}
Note that in taking these upper and lower limits it is enough to consider
 $r \to 0$ through any discrete sequence of $r$ that converges no faster than at a geometric rate. For $d_q^- $ there are convenient alternative forms.
 
\begin{lem}
For $q>1$
\begin{eqnarray} 
d_{q}^{-} & =&
\sup \{s : \limsup_{k\rightarrow \infty}\sum_{\bi \in I_{k}} \svs 
{T_{\bi}}^{1-q}\mu(C_{\bi})^{q} < \infty \}, \label{lgqdim}\\
& =&
\sup \{s :\sum_{k=0}^{\infty} \sum_{\bi \in I_{k}}  
 \svs {T_{\bi}} ^{1-q} 
\mu (C_{\bi})^{q} <\infty \}.\label{seriessum}
\end{eqnarray}
\end{lem}

\noindent{\it Proof.}  
Note  that,  from (\ref{phicont}), 
\be
\phi^{s_1}(T_\bi) \geq   \alpha_{+}^{-k(s-s_1)} \phi^{s}(T_\bi)\label{comp}
\ee
 if
$|\bi|= k$ and $0<s_1<s$. 
Thus if $ \limsup_{k\rightarrow \infty}\sum_{\bi \in I_{k}} \svs 
{T_{\bi}}^{1-q}\mu(C_{\bi})^{q} = M< \infty$ for some $s$, then
$$\sum_{k=0}^{\infty} \sum_{\bi \in I_{k}}  
 \phi^{s_1}(T_{\bi}) ^{1-q} \mu (C_{\bi})^{q}
\leq\sum_{k=0}^{\infty}M\alpha_{+}^{k(s-s_1)(q-1)} <\infty$$
  for all $s_1<s$, and in particular
  $\limsup_{r\rightarrow 0}\sum_{\bi \in J^{s_1}(r)}  \phi^{s_1} 
  (T_{\bi})^{1-q}\mu(C_{\bi})^{q}< \infty$.
  
On the other hand, 
 if $\limsup_{r\rightarrow 0}\sum_{\bi \in J^{s}(r)}  \phi^{s} 
 (T_{\bi})^{1-q}\mu(C_{\bi})^{q}< \infty$ for some $s$, then
 $\sum_{\bi \in J^{s}(r)}  \phi^{s_1} 
 (T_{\bi})^{1-q}\mu(C_{\bi})^{q}\leq  Mr^{a}$
 for some $a>0$, since  (\ref{Jr}) and (\ref{comp}) imply that 
 $\phi^{s_1}(T_\bi) \geq  r^{-c} \phi^{s}(T_\bi)$ if $\bi \in J^{s}(r)$, for some $c>0$.
 If we choose $\rho$ such that $ \alpha_{+}<\rho <1$ then 
 $I = \cup_{k=0}^{\infty} I_k \subseteq \cup_{l=0}^{\infty}J^s(\rho^l)$, so
 $$\sum_{k=0}^{\infty} \sum_{\bi \in I_{k}}  
 \phi^{s_1}(T_{\bi}) ^{1-q} \mu (C_{\bi})^{q}
 \leq  \sum_{l=0}^{\infty} \sum_{\bi \in J^s(\rho^l)}  
 \phi^{s_1}(T_{\bi}) ^{1-q} \mu (C_{\bi})^{q}
 \leq  \sum_{l=0}^{\infty}M \rho^{la} <\infty$$
  for all $s_1<s$, and in particular
 $ \limsup_{k\rightarrow \infty}\sum_{\bi \in I_{k}} \phi^{s_1} 
(T_{\bi})^{1-q}\mu(C_{\bi})^{q} < \infty$.
We conclude that the numbers in (\ref{basicb1}), (\ref{lgqdim}) and (\ref{seriessum}) are equal.
$\Box$

\medskip

There are not, in general, expressions for $d_{q}^{+}$ analogous to (\ref{lgqdim}) and  (\ref{seriessum}). However, we will see in Section 8 that for many measures $\mu$ on $I_{\infty}$, including Bernoulli measures and Gibbs measures,  $d_{q}^{-}= d_{q}^{+}$ so that all these expressions are equal.

It is easy to obtain upper bounds on the $q$-dimensions of the $\mu^{\omega}$ from Proposition \ref{sumest}.

\begin{cor}\label{corub}
Let $\mu$ be a finite Borel measure on $\ii$, let $\mu^{\omega}$ be the measure on $E^\w$ 
defined by {\rm (\ref{mesdef})} and let $q > 1$. Then for all $\w$
$$
\lqd (\mu^\w) \leq  \min\{ d_{q}^{-}(T_{1},\ldots,T_{m};\mu), N\}, 
$$
and
$$
\uqd (\mu^\w) \leq  \min\{d_{q}^{+}(T_{1},\ldots,T_{m};\mu), N\} .
$$
\end{cor}

\noindent{\it Proof.}  
This is immediate from  (\ref{3.2n}) taken in conjunction with the definitions  (\ref{basicb1}) and  (\ref{basicb2}), noting that generalised dimensions never exceed the dimension of the ambient space.
$\Box$
\medskip

\section{The random model}\label{rand} 
\setcounter{equation}{0} 
\setcounter{theo}{0}

Generalised dimensions of measures on (almost) self-affine sets are not everywhere continuous in the defining parameters so we can only hope for generic or almost sure results. More over, as is usually the case, it is harder to get good lower bounds than upper bounds. 

One might hope that generically that one would have $\lqd (\mu^\w)= \min\{d_q^-, N\}$. This was shown to be the case in \cite{Fa5}  for almost all (with respect to translates) self-affine measures with $\|T_i\| < \frac{1}{2}$ for all $i$ in the case $1<q \leq 2$. However, varying the translates did not provide enough randomness in such strictly self-affine constructions to be able to extend this to larger $q$. Here we address this difficulty by working with an almost self-affine model, our ultimate aim being to  show that  $\lqd (\mu^\w)= \min\{d_q^-, N\}$  almost surely for measures on random almost self-affine sets. 

Let $D$ be a bounded region in $ \Esp$. For each $\bi$, let $\w_\bi\in D$ be a random vector distributed according to some Borel probability measure $\P_\bi$ that is absolutely continuous with respect to $N$-dimensional Lebesgue measure. We assume that the $\w_\bi$ are independent  identically distributed random vectors. We let $\P$ denote the product probability measure $\P = \prod_{\bi \in T}  \P_\bi$ on the family of displacements
$\w = \{\w_\bi: \bi \in T\}$. 

In this context, the points $\Pi^\w(\bi) \in \Esp$  given by (\ref{points1})-(\ref{points2}) are now random points whose aggregate form the random set $E^\w$ of  (\ref{2.6}), and the measure  $\mu^{\omega}$ defined by (\ref{mesdef})- (\ref{mesintdef}) is supported by $E^\w$.

The main theorem of \cite{JPS} states that, in this setting, the dimension of $E^\w$ equals the affinity dimension almost surely.

\begin{theo}\label{jorpolsim}
Provided that $\|T_i\|<1 $ for $i=1,\ldots,m$, for almost all $\w$,

$(1)$ $\hdd E^\w = {\mbox{\rm dim}}_{\rm B} (E^{ \w} )= d(T_{1},\ldots, T_{m})$ if $d(T_{1},\ldots, T_{m}) \leq N,$

$(2)$ ${\cal L}^N( E^\w) >0$  if $d(T_{1},\ldots, T_{m}) > N,$

\noindent where  $d(T_{1},\ldots, T_{m})$ is given by $ (\ref{dimexp})$.
\end{theo}

\noindent{\it Proof.} 
This is established  using a potential theoretic method in
 \cite[Theorem 1.5]{JPS}.
$\Box$
\medskip

Our aim now is to obtain an analogue of this result for the generalised $L^q$-dimensions $\lqd(\mu^{\w})$ of the measures 
$\mu^{\omega}$ for $q>1$. The upper estimate was addressed in Section 4. For the lower estimate we proceed  in two stages. We first obtain an upper bound for $\E \int \mu^{\omega}  (B(x,r))^{q-1}  d\mu^{\omega}(x)$, the expectation of the quantity that occurs in the definition of the generalised dimensions  (\ref{intdefl}), in terms of a  `multienergy integral' (\ref{probest2}). We then use an induction on trees to show that this  integral is bounded if $\sum_{|\bi| =k} \ps(T_\bi)^{1-q} \mu(C_{\bi})^q \to 0$ geometrically as $k \to \infty$.

For the first stage, we recall that the inverse singular values $\svs {T_{\bi\wedge \bj}}^{-1}$, which depend on the join $\bi\wedge \bj$ of $\bi,\bj \in I_\infty$,
play an important r\^{o}le in estimates involving dimensions of self-affine sets.
Multienergy kernels may be regarded as a generalisation of
such expressions to  several points of $I_\infty$. The {\em join set} of $\bi_1,\ldots,\bi_n \in I_\infty$ is the set of  {\em join points} 
\be
J \equiv  \bigwedge(\bi_1,\ldots,\bi_n) = \{\bi_p\wedge\bi_q: p\neq q\}\label{joinset}
 \ee 
with repetitions counted by multiplicity in a natural way. The multiplicity of  
$\bv \in J$ is $r-1$, where $r$ is the greatest integer such that there are distinct  $\bi_{i_1},\ldots,\bi_{i_r}$ with $\bi_{i_p}\wedge \bi_{i_q} = \bv$ for all $1 \leq p<q \leq r$; this ensures that $J$ contains exactly $n-1$ points including repetitions. In the simplest case where $m=2$, every vertex of a join set has multiplicity $1$. It is natural to think of the join points as vertices of the $m$-ary tree $I$ where the paths from $\emptyset$ to the $\bi_j$ meet, see Figure 1.

\begin{figure}
	\centering
	\includegraphics[width=80mm]{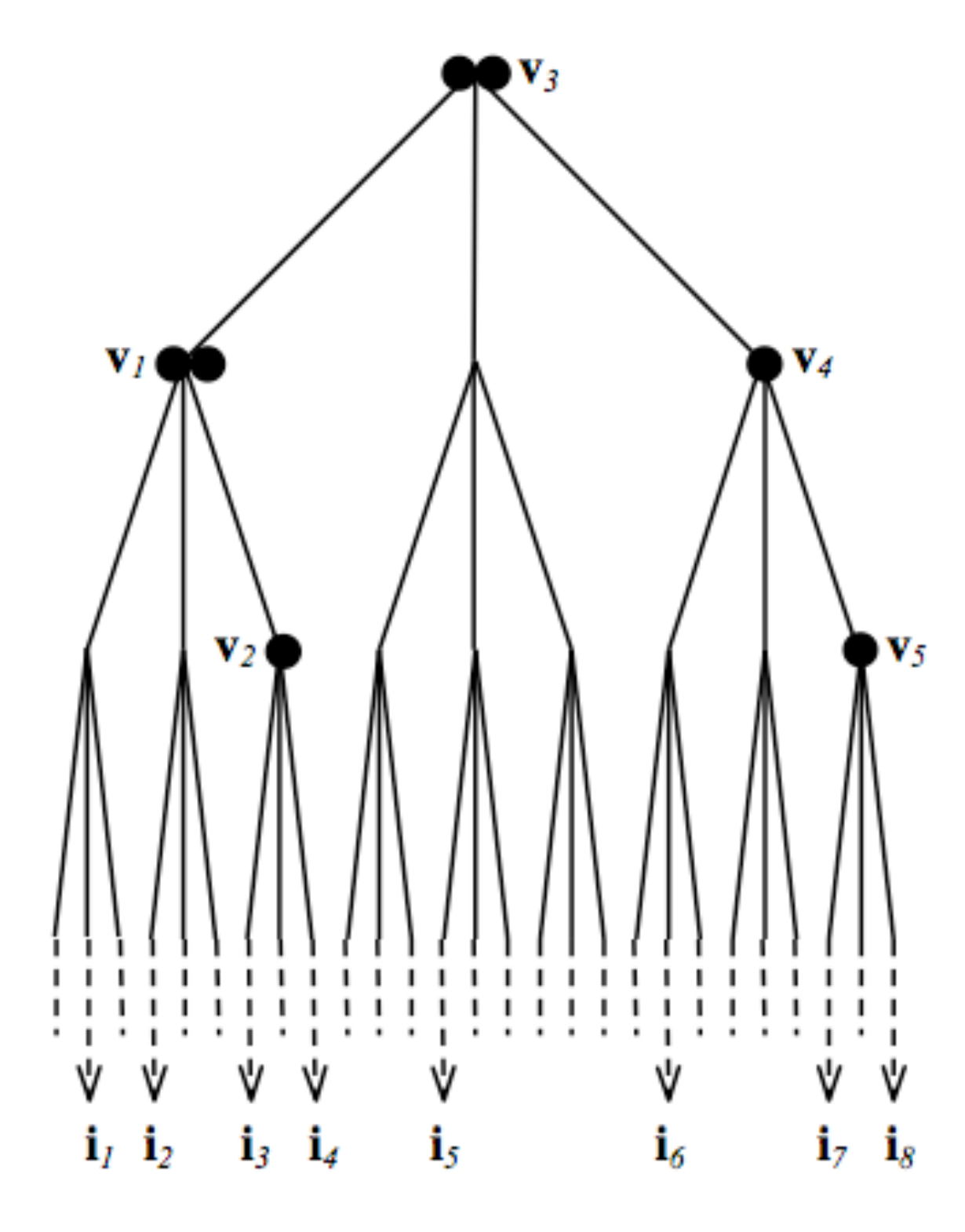}
	\caption{The join set $ \{\bv_1, \bv_1,\bv_2,\bv_3,\bv_3,\bv_4,\bv_5\}$ of 
$ \{\bi_1, \bi_2,\bi_3,\bi_4,\bi_5,\bi_6,\bi_7,\bi_8\}$}
	\label{figure1.pdf}
\end{figure}

We define  multienergy kernels by forming  products of the singular value functions at the vertices of join sets. For 
$\bi_1,\ldots,\bi_n \in I_\infty$ let
\be
 \phi^s(\bi_1,\ldots,\bi_n) =  \phi^s(T_{\bv_1})  \phi^s(T_{\bv_2})\cdots \phi^s(T_{\bv_{n-1}})
 \mbox{ where } \{\bv_1,\ldots,\bv_{n-1}\} =   \bigwedge(\bi_1,\ldots,\bi_n) .\label{multen}
\ee
We will consider  multienergy integrals of the form
$$\int\cdots\int  \phi^s(\bi_ 1,\ldots,\bi_n)^{-1} d\mu(\bi_1)\ldots d\mu(\bi_n)$$
which provide bounds for the expectation $\E \int \mu^{\omega}  (B(x,r))^{q-1}  d\mu^{\omega}(x)$.

The second stage involves showing that,  for suitable $s$, certain multienergy integrals are finite, implying  that  $\int\mu^{\omega}  (B(x,r))^{q-1}  d\mu^{\omega}(x)<\infty$ almost surely, to give a lower bound for 
$\lqd (\mu^w)$.

These two stages are executed in the next two sections.

\section{Probabilistic estimates} 
\setcounter{equation}{0} 
\setcounter{theo}{0}

The aim of this section is to  bound the expectation of $ \int \mu^{\omega}  (B(x,r))^{q-1}  d\mu^{\omega}(x)$, that is the integral which occurs in the definition of the generalised dimensions  (\ref{intdef}), in terms of a multienergy integral.

Let $\E$ denote expectation. Given $J \subseteq I$ we write ${\cal F} = \sigma\{\w_\bi : \bi \in J\}$  for the sigma-field generated by the random displacements $\w_\bi \in J$ and  write $\E(Z \mid {\cal F})$ for the expectation of a random variable $Z$ conditional on  ${\cal F}$; intuitively this is the expectation of $Z$ given all the displacements $\{\w_\bi : \bi \in J\}$. 

Note that the constant $c$ may differ in each of the following lemmas. The first lemma is a `transversality' property  of a form often encountered in work on self-affine sets.

\begin{lem}\label{prolem}
Let $0<s \leq N$ with $s$ not an integer. Then there exists $c>0$ such that
\begin{equation*}
\E( |\xw(\bu)-\xw(\bv) |^{-s}\mid {\cal F}) \leq c \phi^s(T_{\bu\wedge \bv})^{-1}
 \end{equation*}
for all $\bu,\bv \in I$, where
${\cal F} = \sigma\{\w_\bi : \bi \in J\}$ for any subset $J$ of $T$ such that
 $ \bv|_{k+1}, \bv|_{k+2},\ldots \in J$  and  $ \bu|_{k+2}, \bu|_{k+3},\ldots \in J$ but   
 $ \bu|_{k+1}\notin J$, where $|\bu \wedge \bv | = k$.
\end{lem}

\noindent{\it Proof.} 
From (\ref{points2}), for each $\bu,\bv$,
\begin{align*}
 \xw(\bu)-\xw(\bv) 
= T_{\bu\wedge \bv}\big((&\w_{\bu|_{k+1}}+ T_{\bu|_{k+1}}\w_{\bu|_{k+2}} 
+ T_{\bu|_{k+1}}T_{\bu|_{k+2}}\w_{\bu|_{k+3}} + \cdots)\\
& -(\w_{\bv|_{k+1}}+ T_{\bv|_{k+1}}\w_{\bv|_{k+2}} 
+ T_{\bv|_{k+1}}T_{\bv|_{k+2}}\w_{\bv|_{k+3}} + \cdots)\big)\\
= T_{\bu\wedge \bv}\big((&\w_{\bu|_{k+1}}+ x(\w)\big)
\end{align*}
where $x(\w)$ is ${\cal F}$-measurable. Thus
\begin{align*}
\E( |\xw(\bu)-\xw(\bv) |^{-s}\mid {\cal F}) 
= & \int\frac{d\P( \w_{\bu|_{k+1}})}{|T_{\bu\wedge \bv}\big((\w_{\bu|_{k+1}}+ x(\w)\big) |^{s}}\\
\leq & c \phi^s(T_{\bu\wedge \bv})^{-1},
\end{align*}
where the integral may be estimated just as in \cite[Lemma 3.1]{Fa1} or \cite[Lemmas 4.5, 5.2]{JPS}.
$\Box$
\medskip

We next use a sequence of conditional expectations to extend Lemma \ref{prolem} from $2$ to $n+1$ points of $I_\infty$. 

\begin{lem}\label{prolem2}
For all $0<s \leq N$ with $s$ not an integer,  there exist numbers $c>0$ and $r_0>0$ such that for all  $\bi_1,\ldots,\bi_n,\bj \in I_{\infty}$ and $0<r \leq r_0$,
\begin{equation}
\P\big\{|\xw(\bi_1)-\xw(\bj)|  \leq r,\ldots,|\xw(\bi_n)-\xw(\bj)| \leq r\big\}
\leq c r^{sn}\phi^s(\bi_1,\ldots,\bi_n,\bj)^{-1}.
\label{probdists}
\end{equation}
\end{lem}
 
\noindent{\it Proof.} 
We may renumber the points $\bi_1,\ldots,\bi_n$ in such a manner that  
$\{\bi_1\wedge \bi_2, \bi_2\wedge \bi_3,\ldots, \bi_{n-1}\wedge \bi_n, \bi_n\wedge \bj\} $ are precisely the points of the join set $\bigwedge(\bi_1,\ldots,\bi_n)$, including any repeated points. (One way to achieve this renumbering is to transform the tree $I$ by an automorphism fixing the root $\emptyset$ in such a way that $\bj$ is the  \textquoteleft extreme right' point of the tree and renumber the $\bi_k$ from left to right.) Note that this renumbering does not affect the value of $\phi^s(\bi_1,\ldots,\bi_n,\bj)$. Thus
\begin{align}
\P\big\{ &|\xw(\bi_1) -\xw(\bj)|  \leq r,\ldots,|\xw(\bi_n)-\xw(\bj)| \leq r\big\}\nonumber\\
& \leq \P\big\{|\xw(\bi_1)-\xw(\bi_2)|   \leq 2r,\ldots,|\xw(\bi_{n-1})-\xw(\bi_n)| \leq 2r,
|\xw(\bi_n)-\xw(\bj)| \leq 2r\big\}\nonumber\\
& \leq 2^nr^{sn}
\E\big( |\xw(\bi_1)-\xw(\bi_2)|^{-s} \cdots |\xw(\bi_{n-1})-\xw(\bi_n)|^{-s}
|\xw(\bi_n)-\xw(\bj)|^{-s}\big). \label{expint}
\end{align}
We estimate this expectation through a tower of conditional expectations.
Define a sequence of sigma-fields ${\cal F}_1 \supset {\cal F}_2 
\supset \cdots \supset {\cal F}_n$ by ${\cal F}_l = \sigma\{\w_\bi : \bi \neq \bi_1 |_{k_1+1},\ldots, \bi_l |_{k_l+1}\}$
where $k_l = |\bi_l\wedge \bi_{l+1}| \, (1 \leq l \leq n-1)$ and   $k_n = |\bi_n \wedge \bj|$.

For brevity of notation, write $Z^{\w}_l = |\xw(\bi_l)-\xw(\bi_{l+1})|^{-s}\, (1 \leq l \leq n-1)$ with 
$Z^{\w}_n = |\xw(\bi_n)-\xw(\bj)|^{-s}$, so that $Z^{\w}_{l+1},\ldots, Z^{\w}_n$ are all ${\cal F}_{l}$-measurable for $l =1,\ldots,n-1$. Using the tower property for conditional expectation, that  $Z^{\w}_2\ldots Z^{\w}_n$ is 
${\cal F}_1$-measurable, and then applying Lemma \ref{prolem},
\begin{align*}
\E\big(Z^{\w}_1\ldots Z^{\w}_n \big|{\cal F}_n \big)
&= \E\Big(\E\big(Z^{\w}_1\ldots Z^{\w}_n\big|{\cal F}_1\big)\big|{\cal F}_n \Big)\\
&= \E\Big(\E\big(Z^{\w}_1\big|{\cal F}_1\big)Z^{\w}_2\ldots Z^{\w}_n\big|{\cal F}_n \Big)\\
&\leq \E\big(c \phi^s(T_{\bi_1\wedge \bi_2})^{-1} Z^{\w}_2\ldots Z^{\w}_n\big|{\cal F}_n \big)\\
&= c \phi^s(T_{\bi_1\wedge \bi_2})^{-1} \E\big(Z^{\w}_2\ldots Z^{\w}_n\big|{\cal F}_n \big).
\end{align*}
Repeating this argument $n-1$  times, we obtain that
$$\E\big(Z^{\w}_1\ldots Z^{\w}_n \big|{\cal F}_n \big)
\leq c^n  \phi^s(T_{\bi_1\wedge \bi_2})^{-1}\ldots  \phi^s(T_{\bi_{n-1}\wedge \bi_n})^{-1}
 \phi^s(T_{\bi_n\wedge \bj})^{-1} = c^n   \phi^s(\bi_1,\ldots,\bi_n,\bj)^{-1},$$
giving the bound for the unconditional expectation
\begin{align*}
\E\big( |\xw(\bi_1)-\xw(\bi_2)|^{-s} &\cdots|\xw(\bi_{n-1})-\xw(\bi_n)|^{-s}
|\xw(\bi_n)-\xw(\bj)|^{-s} \big)\\
&=
\E\big(Z^{\w}_1\ldots Z^{\w}_n) \leq  c^n  \phi^s(\bi_1,\ldots,\bi_n,\bj)^{-1}.
\end{align*}
Combining this with (\ref{expint}) gives (\ref{probdists}). 
$\Box$

\medskip

We now integrate (\ref{probdists}) over the $\bi_l$.

\begin{lem}
For all $0<s \leq N$ with $s$ not an integer,  there exist numbers $c>0$ and $r_0>0$ such that for all $\bj \in I_{\infty}$ and  $0<r\leq r_0$,
\begin{equation}
\E\big( \mu^{\omega}(B(\xw(\bj),r))^n\big)  \leq c r^{sn} 
\int \cdots \int 
 \phi^s(\bi_1,\ldots,\bi_n,\bj)^{-1} d\mu(\bi_1)\ldots d\mu(\bi_n).\label{probest1}
\end{equation}
 \end{lem} 
 
\noindent{\it Proof.} 
Using Fubini's theorem, 
\begin{align*}
\E\big( \mu^{\omega}&(B(\xw(\bj), r))^n\big)\\
& = \E\big(\mu\{\bi :  |\xw(\bi)-\xw(\bj) | \leq r\}^n \big)\\ 
&= \big( \P\times \mu \times \cdots\times \mu \big) \big\{(\omega,\bi_1,\ldots,\bi_n):
 |\xw(\bi_1)-\xw(\bj)|  \leq r,\ldots,|\xw(\bi_n)-\xw(\bj)| \leq r\big\}\\ 
 &= \int\cdots\int \P\big\{
 |\xw(\bi_1)-\xw(\bj)|  \leq r,\ldots,|\xw(\bi_n)-\xw(\bj)| \leq r\big\} d\mu(\bi_1)\ldots d\mu(\bi_n)\\
 &\leq c r^{sn} 
\int \cdots \int 
 \phi^s(\bi_1,\ldots,\bi_n,\bj)^{-1} d\mu(\bi_1)\ldots d\mu(\bi_n),
\end{align*}
by Lemma \ref{prolem2}.
 $\Box$

\medskip

Finally, integration of an appropriate power of (\ref{probest1}) with respect to  $\bj$ gives the main estimate (note the simpler form of (\ref{probest2}) if $q=n+1$ is an integer).

\begin{prop}\label{propexp}
Let $n \geq 1$ and $1 < q \leq n+1$. Then for all $0<s \leq N$ with $s$ not an integer, 
there exist numbers $c>0$ and $r_0>0$ 
such that for all $0<r\leq r_0$, 
\begin{align}
\E \int \mu^{\omega} & (B(x,r))^{q-1}  d\mu^{\omega}(x) \nonumber\\
&\leq c r^{s(q-1)} 
 \int\bigg[
 \int \cdots \int 
 \phi^s(\bi_1,\ldots,\bi_n,\bj)^{-1} d\mu(\bi_1)\ldots d\mu(\bi_n)\bigg]^{(q-1)/n} d\mu(\bj)\label{probest2}
\end{align}
\end{prop} 
 
\noindent{\it Proof.} 
Since $n/(q-1) \geq 1$, Jensen's inequality and (\ref{probest1}) give
\begin{align*}
\E \int \mu^{\omega}(B(x, & r))^{q-1}  d\mu^{\omega}(x) \\
& = \E \int \mu^{\omega}  (B(\xw(\bj),r))^{q-1}  d\mu(\bj) \\
& \leq  \int\Big[\E\big( \mu^{\omega}  (B(\xw(\bj),r))^{n}\big)  \Big]^{(q-1)/n}d\mu(\bj) \\
&\leq c r^{s(q-1)} 
 \int\bigg[
 \int \cdots \int 
 \phi^s(\bi_1,\ldots,\bi_n,\bj)^{-1} d\mu(\bi_1)\ldots d\mu(\bi_n)\bigg]^{(q-1)/n} d\mu(\bj).
\end{align*}
$\Box$

\section{Integral estimates}\label{integralest}
\setcounter{equation}{0}
\setcounter{theo}{0}

This section is devoted to estimating the integral (\ref{probest2}). To do this we identify the code space 
$I$ with the vertices of a rooted $m$-ary
tree with root $\emptyset$, in the obvious way. Thus the edges of the tree  join  each $\bi \in I$  to its $m$  `children'   $\bi 1,\ldots,\bi m$. To  estimate (\ref{probest2}) we will split the domain of integration into subdomains consisting of  $n$-tuples $(\bi_1 ,\ldots,\bi_n)$ whose join sets lie in certain families of automorphisms of the tree $I$. We will use induction over classes of join sets to estimate the integrals over each such domain, with H\"{o}lder's inequality playing a very natural r\^{o}le at each step.

We require a little terminology.
A {\it join set} $J$ with {\it root} $\bw \in I$ consists of a family of vertices $\{\bv_1,\ldots,\bv_r\}$ of $I$, with repetitions allowed,  such that $\bv_i \succeq \bw$ for all $i$ and with the property that
$\bv_i\wedge \bv_j \in J$ for all $\bv_i,\bv_j \in J$. The root may or may not be a vertex of the join set. The number $r+1$ is called the {\it spread} of the join set. The  {\it multiplicity} of a given vertex is the number of times it occurs in $J$.

Join sets occur naturally in connection with the integrals (\ref{probest2}): given $\bi_1,\ldots,\bi_n \in I_\infty$ then $J = \bigwedge(\bi_1,\ldots,\bi_n)$ is a join set of spread $n$, and the multiplicity of a vertex 
$\bv \in J$ is $r-1$ where $r$ is the greatest integer such that there are distinct  $\bi_{i_1},\ldots,\bi_{i_r}$ with $\bi_{i_p}\wedge \bi_{i_q} = \bv$ for all $1 \leq p<q \leq r$. (In a binary tree every vertex of a join set has multiplicity 1.)

A {\it join class} ${\cal J}$ with {\it root} $\bw \in I$ is an equivalence class of join sets all with root $\bw$, two such join sets $J$ and $J'$ being equivalent if there is an automorphism of the rooted subtree of $I$ with root $\bw$ that maps $J$ onto $J'$ (preserving multiplicities). The {\it spread} of  a join class ${\cal J}$ is the common value of the spreads of all  $J \in {\cal J}$.

The {\it level} of a vertex $\bv \in I$ is just $|\bv|$.  Thus the  {\it set of levels} $L(J)$ of a join set $J=
\{\bv_1,\ldots,\bv_r\}$ is  $\{|\bv_1|,\ldots,|\bv_r|\}$, allowing repetitions, and the  {\it set of levels} $L({\cal J})$ of a join class is the common set of levels of the join sets in the class.

Note that in   (\ref{intest}) and below, the product is over the set of levels in a join class. The symbol $[n-1]$ above the product sign merely indicates that there are  $n-1$ terms in this product;  this convention is helpful when keeping track of terms through the proofs.

\begin{prop}\label{integerest}
Let $q>1$ and $n\geq 2$ be such that $q\geq n$.  Let $\J$ be a join class with root $\bv$ and spread $n$.  Then
\begin{align}
\int_{\bigwedge(\bi_1,\ldots,\bi_n) \in \J}  & \phi^s(\bi_1,\ldots,\bi_n)^{-1} d\mu(\bi_1)\ldots d\mu(\bi_n)\nonumber  \\ 
& \leq \mu(C_\bv)^{(q-n)/(q-1)} \prod_{l \in L(\J)}^{[n-1]}
 \Big(\sum_{|\bu|= l, \bu \succeq \bv} \ps(T_\bu)^{1-q} \mu(C_{\bu})^q\Big)^{1/(q-1)}. \label{intest}
\end{align}
\end{prop}

\noindent{\it Proof.} 
We proceed by induction on the number of {\it distinct} vertices of $\J$. To start the inductive process, suppose that the join sets in $\J$ consist of a single vertex $\bv$ of multiplicity $n-1$ for some $n \geq 2$. First assume that  $\bv$ is itself the root of the join sets of $\J$. Then
\begin{align}
\int_{\bigwedge(\bi_1,\ldots,\bi_n) \in \J}  \phi^s(\bi_1,\ldots,\bi_n)^{-1} &d\mu(\bi_1)\ldots d\mu(\bi_n)
 = \int_{\bigwedge(\bi_1,\ldots,\bi_n) \in \J}   \phi^s(T_\bv)^{-(n-1)} d\mu(\bi_1)\ldots d\mu(\bi_n)\nonumber\\
&\leq    \phi^s(T_\bv)^{-(n-1)} \mu(C_\bv)^n\nonumber\\
&=  \mu(C_\bv)^{(q-n)/(q-1)} \big( \phi^s(T_\bv)^{1-q} \mu(C_\bv)^q\big)^{(n-1)/(q-1)},\label{firstest}
\end{align}
which is  (\ref{intest}), noting that $L({\cal J})$ has just  one level $ |\bv|$ which is of multiplicity $n-1$, with the sums in each multiplicand of (\ref{intest})  having just one term each. 

If, now, the join class $\J$ has root $\bw$ and contains join sets $J$ consisting of a single vertex $\bv$, distinct from $\bw$, of multiplicity $n-1$ at level $l$, we may sum (\ref{firstest}) over  $\bv$ such that $\bv \succeq \bw$ and 
$|\bv|=l$
to get, using H\"{o}lder's inequality,
\begin{align*}
\int_{\bigwedge(\bi_1,\ldots,\bi_n) \in \J}  & \phi^s(\bi_1,\ldots,\bi_n)^{-1} d\mu(\bi_1)\ldots d\mu(\bi_n)\\
&\leq \sum_{|\bv| = l, \bv \succeq \bw}  \mu(C_\bv)^{(q-n)/(q-1)} \big( \phi^s(T_\bv)^{1-q} \mu(C_\bv)^q\big)^{(n-1)/(q-1)}\\
&\leq \Big(\sum_{|\bv| = l, \bv \succeq \bw}  \mu(C_\bv)\Big)^{(q-n)/(q-1)} \Big( \sum_{|\bv| = l, \bv \succeq \bw}\phi^s(T_\bv)^{1-q} \mu(C_\bv)^q\Big)^{(n-1)/(q-1)}\\
&= \mu(C_\bw)^{(q-n)/(q-1)} \Big( \sum_{|\bv| = l, \bv \succeq \bw}\phi^s(T_\bv)^{1-q} \mu(C_\bv)^q\Big)^{(n-1)/(q-1)},
 \end{align*}
which is  (\ref{intest}) for join classes with a single vertex of any multiplicity.

Now assume inductively that  (\ref{intest}) holds for all join sets with fewer than $k$ distinct vertices for some $k\geq 2$. Let $\J$ be a join class with root $\bv$ and $k\geq 2$ distinct vertices and spread $n$ where $n \leq q$. Again, first consider the case where that the root $\bv$  belongs to the join sets  in $\J$ as the \textquoteleft  top' vertex. In each join set $J \in \J$ there is a (possibly empty) set of $r \geq 0$ vertices $\{\bv_1,\ldots,\bv_r\}$ in $J$ distinct from and  \textquoteleft immediately below' $\bv$, that is with the path joining $\bv_i$ to $ \bv$  in the tree $I$ containing no other vertices of $J$.  For a given class $ \J$ these sets of vertices (with multiplicity) are equivalent under automorphisms of the tree that fix the root $\bv$.

For each $i$, the join set  $J \in \J$ induces a join set that we denote by $J_i$ with root $\bv$ and vertices 
$\{ \bu \in J: \bu \succeq \bv_i\}$, that is the vertices of $\J$ below and including $\bv_i$. These join sets are equivalent under automorphisms of the tree that fix the root $\bv$ and we write $\J_i = \{J_i : J \in {\cal J} \}$ for this equivalence class of join sets, which has spread $n_i\geq 2$, say,  and set of levels $L_i$ (counted with repetitions).

Let 
\be
n_1 +\cdots +n_r +t = n, \label{weightsum}
\ee
where $t\geq 0$. To integrate over $\{\bi_1,\ldots,\bi_n\}$ such that $\bigwedge(\bi_1,\ldots,\bi_n)\in  \J$ we decompose the integral so that for each $J \in  \J$,  for each $i \, (1 \leq i \leq r)$, the $n_i $ integration variables, 
$\{\bi_1^i,\ldots,\bi_{n_i}^i\}$ say,
are such that $\bigwedge(\bi_1^i,\ldots,\bi_{n_i}^i) = J_i$, and $t$  of them, $\{\bi_1^0,\ldots,\bi_{t}^0\}$, such that $\bi_l^0\wedge \bi_r = \bv$ for all $ \bi_r \neq \bi_l^0$. Thus, noting that the multiplicity of $\bv$ is $r+t-1$,
\begin{align*}
I  &\equiv \int_{\bigwedge(\bi_1,\ldots,\bi_n) \in \J}\phi^s(\bi_1,\ldots,\bi_n)^{-1} d\mu(\bi_1)\ldots d\mu(\bi_n)\\
& 
\leq \phi^s(T_\bv)^{-(r+t-1)}
\mu(C_\bv)^t 
\int_{\bigwedge(\bi_1^1,\ldots,\bi_{n_1}^1) \in \J_1}  
 \phi^s(\bi_1^1,\ldots,\bi_{n_1}^1)^{-1} d\mu(\bi_1^1)\ldots d\mu(\bi_{n_1}^1)\\
& \quad\quad \times \cdots \times
\int_{\bigwedge(\bi_1^r,\ldots,\bi_{n_r}^r) \in \J_r }  \phi^s(\bi_1^r,\ldots,\bi_{n_r}^r)^{-1} d\mu(\bi_1^r)\ldots d\mu(\bi_{n_r}^r)\\
& 
\leq \phi^s(T_\bv)^{1-r-t} \mu(C_\bv)^t  
\mu(C_{\bv})^{(q-n_1)/(q-1)} \prod_{l \in L_1}^{[n_1-1]} 
 \Big(\sum_{|\bu| =l, \bu \succeq \bv} \ps(T_\bu)^{1-q} \mu(C_{\bu})^q\Big)^{1/(q-1)} \\
& \quad\quad  \times \cdots \times
\mu(C_{\bv})^{(q-n_r)/(q-1)} \prod_{l \in L_r}^{[n_r-1]}
 \Big(\sum_{|\bu| =l, \bu \succeq \bv} \ps(T_\bu)^{1-q} \mu(C_{\bu})^q\Big)^{1/(q-1)},
\end{align*}
where we have applied  the inductive assumption (\ref{intest}) to join sets in $\J_i$ (which all have root $\bv$) for each $i$. 

Combining terms,
\begin{align*}
I  &\leq \mu(C_{\bv})^{(q-n_1- \cdots- n_r -t)/(q-1)} \big(\phi^s(T_\bv)^{1-q} \mu(C_\bv)^q\big)^{(r+t-1)/(q-1)}\\
& \qquad\qquad  \qquad\qquad \qquad\times 
\prod_{l \in L_1\cup \cdots \cup  L_r}^{[n_1+ \cdots+ n_r-r]} 
 \Big(\sum_{|\bu| =l, \bu \succeq \bv} \ps(T_\bu)^{1-q} \mu(C_{\bu})^q\Big)^{1/(q-1)} \\
&\leq \mu(C_{\bv})^{(q-n)/(q-1)} 
\prod_{l \in L}^{[n-1]} 
 \Big(\sum_{|\bu| =l, \bu \succeq \bv} \ps(T_\bu)^{1-q} \mu(C_{\bu})^q\Big)^{1/(q-1)},
 \end{align*}
where $L$ is the complete set of levels of  ${\cal J}$, where we have used (\ref{weightsum}), and incorporated the terms $\phi^s(T_\bv)^{1-q} \mu(C_\bv)^q$ (taken as a sum over the single vertex $\bv$) in the main product with multiplicity $(r+t-1)$.  This is  (\ref{intest}) in the case where the root of the join sets in ${\cal J}$ belongs to the join sets.

Finally, if the root $\bw$ is not a vertex of the join sets in ${\cal J}$, then summing (\ref{intest}) over join sets with top vertex $\bv \preceq \bw$ with $|\bv| = l'$ and using H\"{o}lder's inequality,

\begin{align*}
\int_{\bigwedge(\bi_1,\ldots,\bi_n) \in \J}  & \phi^s(\bi_ 1,\ldots,\bi_n)^{-1} d\mu(\bi_1)\ldots d\mu(\bi_n)\nonumber 
\\ 
& \leq\sum_{|\bv| = l', \bv \succeq \bw}   \mu(C_\bv)^{(q-n)/(q-1)} \prod_{l \in L(\J)}^{[n-1]}
 \Big(\sum_{|\bu| =l, \bu \succeq \bv} \ps(T_\bu)^{1-q} \mu(C_{\bu})^q\Big)^{1/(q-1)}\\
& \leq \Big(\sum_{|\bv| = l', \bv \succeq \bw}   \mu(C_\bv)\Big)^{(q-n)/(q-1)} \prod_{l \in L(\J)}^{[n-1]}
 \Big( \sum_{|\bv| = l', \bv \succeq \bw} \sum_{|\bu| = l, \bu \succeq \bv} \ps(T_\bu)^{1-q} \mu(C_{\bu})^q\Big)^{1/(q-1)}\nonumber\\
 & = \mu(C_\bw)^{(q-n)/(q-1)} \prod_{l \in L(\J)}^{[n-1]}
 \Big(  \sum_{|\bu| = l, \bu \succeq \bw} \ps(T_\bu)^{1-q} \mu(C_{\bu})^q\Big)^{1/(q-1)},
  \end{align*}
which completes the induction and the proof.
$\Box$

\medskip

Proposition \ref{integerest} would be adequate for our needs for the cases when $q$ is an integer. However for  non-integral  $q>1$ we need a generalisation where one of the variables of integration is distinguished. 
Whilst the proof of Proposition \ref{fracest} again uses induction on join sets and H\"{o}lder's inequality, the details are more intricate than in Proposition \ref{integerest}, and indeed depends on Proposition \ref{integerest} at several points.

Let $0 \leq k_1 < k_2 < \ldots <k_p$ be  levels, where $ 1\leq p \leq n$. For $\bj \in I_{\infty}$ write $\bj_r =\bj|_{k_r}$. For each $r=1,\ldots,p$, let $\J_r$ be a given join class with root at level $k_r$ and spread $m_r$, see Figure 2. Write $\J_r(\bj_r)$ for the corresponding join class with the root    
mapped to $\bj_r$, that is so that the tree automorphisms of $I$ that map $\bj_r$ to $\bj_r'$ map the join sets in $\J_r(\bj_r)$ onto those in $\J_r(\bj_r ')$ in a bijective manner.

We need to include the cases where $\J_r(\bj_r)$ is a join class with root $\bj_r$ and of spread $1$. In this case we  interpret integration over
$\bigwedge(\bi_{l}) \in \J_{r}(\bj_r)$ as integration over all $\bi_l$ such that  $\bi_l\wedge \bj = \bj_r$, and we take $\phi^s(\bi_{l})= 1$ where this occurs in the next proof.

\begin{prop}\label{fracest}
Let $q>1$ and let $n$ be an integer with $n\geq q-1$. With notation as above, for each $r=1,\ldots,p$ let $\J_r$ be a join class with root at level $k_r$ and spread $m_r  \geq 1$, with $m _1+\cdots + m_p = n$.   
Then
\begin{align}
\int_{\bj \in I}\bigg[\int_{\bigwedge(\bi_{1},\ldots,\bi_{m_{1}}) \in \J_{1}(\bj_1)} \cdots
\int_{\bigwedge(\bi_{n-m_p +1},\ldots,\bi_{n}) \in \J_{p}(\bj_p)} &
 \phi^s(\bi_1,\ldots,\bi_n,\bj)^{-1} d\mu(\bi_1)\ldots d\mu(\bi_n)\bigg]^{(q-1)/n} d\mu(\bj) \nonumber  \\ 
& \leq 
\prod_{l \in L}^{[n]}
 \Big(\sum_{|\bu| =l} \ps(T_\bu)^{1-q} \mu(C_{\bu})^q\Big)^{1/ n}, \label{qintest}
\end{align}
where $ L$ denotes the aggregate set of levels of $\{L(\J_1),\ldots, L(\J_p),k_1,\ldots,k_p\}$.
\end{prop}

\noindent{\it Proof.} 
We proceed by induction on $r $, starting with $r= p$ and working up to $r=1$, taking as the inductive hypothesis:

\noindent For all $\bj_r \in I_{ k_r}$, 
\begin{align}
\int_{\bj \succeq \bj_r}\bigg[
& \int_{\bigwedge(\bi^{r}_{1},\ldots,\bi^{r}_{m_{r}}) \in \J_{r}(\bj_{r})}\cdots
\int_{\bigwedge(\bi^{p}_{1},\ldots,\bi^{p}_{m_{p}}) \in \J_{p}(\bj_p)}
 \phi^s(\bi^{r}_{1},\ldots,\bi^{r}_{m_{r}})^{-1} \phi^s(T_{\bj_{r}})^{-1}\times \cdots  \nonumber \\
&\qquad \times \phi^s(\bi^{p}_{1},\ldots,\bi^{p}_{m_{p}})^{-1} \phi^s(T_{\bj_{p}})^{-1}
 d\mu(\bi^{r}_{1})\ldots d\mu(\bi^{r}_{m_r})\cdots 
 d\mu(\bi^{p}_{1})\ldots d\mu(\bi^{p}_{m_p})\bigg]^{(q-1)/n} d\mu(\bj) \nonumber  \\ 
& 
\leq \mu(C_{\bj_r}) ^{(n - n_r )/n}
\prod_{l \in L_r}^{[ n_r]}
 \Big(\sum_{|\bu| =l, \bu \succeq \bj_r } \ps(T_\bu)^{1-q} \mu(C_{\bu})^q\Big)^{1/n},\label{indhop}
\end{align}
where $n_r = m_r+\cdots+m_p$ and  $ L_r$ denotes the set of levels of $\{L(\J_r),\ldots, L(\J_p),k_r,\ldots,k_p\}$ counted by multiplicity (so that  $ L_r$ consists of $m_r+\cdots+m_p = n_r-1$ levels).

\begin{figure}
	\centering
	\includegraphics[width=100mm]{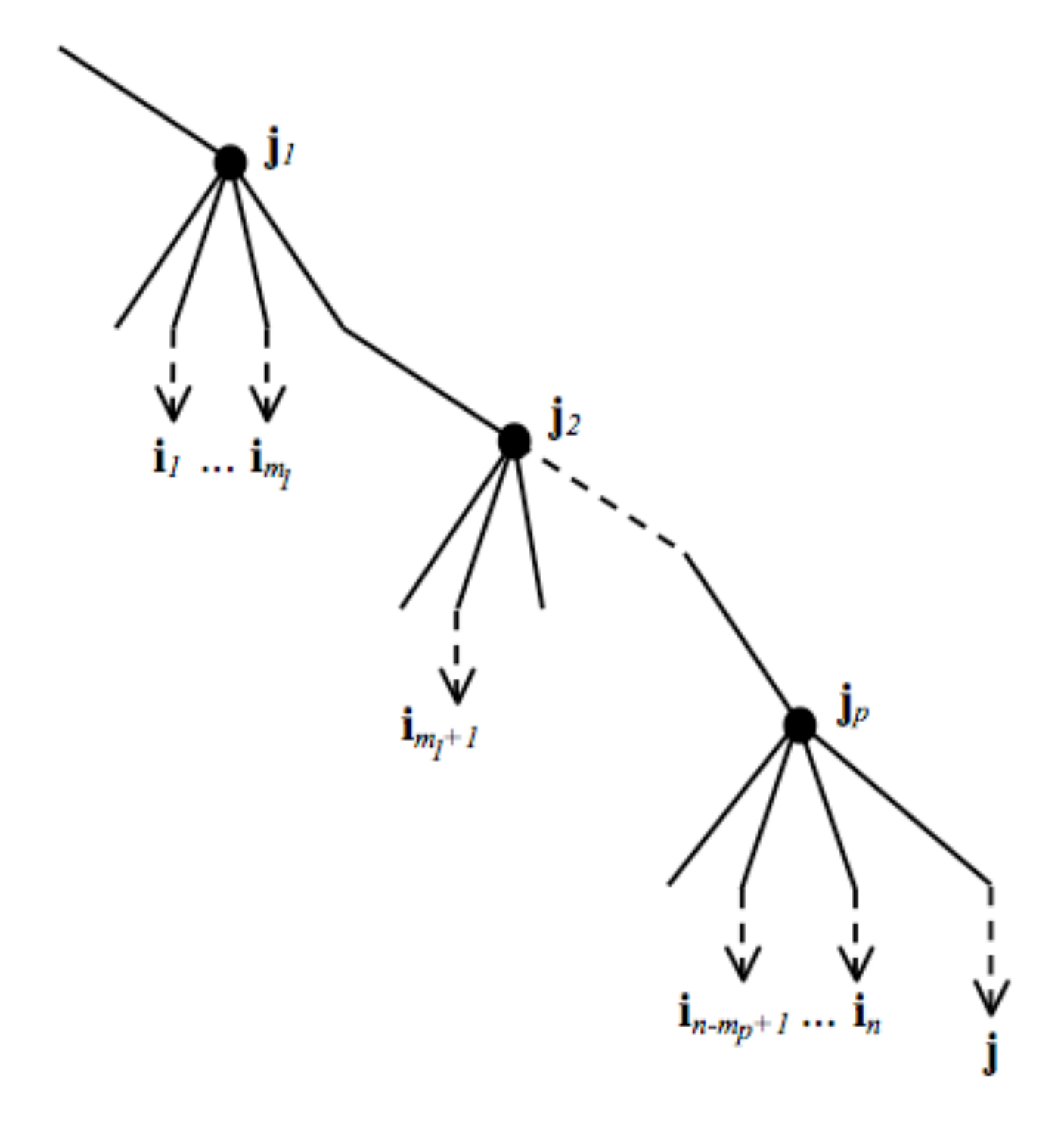}
	\caption{Part of the tree $I$ indicating the notation of Proposition 7.2}
	\label{figure2}
\end{figure}

To start the induction, we apply Proposition \ref{integerest} to get 

\begin{align}
\int_{\bj \succeq \bj_p}&\bigg[
\int_{\bigwedge(\bi^{p}_{1},\ldots,\bi^{p}_{m_{p}}) \in \J_{p}(\bj_p)}
 \phi^s(\bi^{p}_{1},\ldots,\bi^{p}_{m_{p}})^{-1} \phi^s(T_{\bj_{p}})^{-1}
 d\mu(\bi^{p}_{1})\ldots d\mu(\bi^{p}_{m_p})\bigg]^{(q-1)/n} d\mu(\bj)  \label{exp1}\\ 
&
\leq  \int_{\bj \succeq \bj_p}\bigg[
\mu(C_{\bj_p}) ^{(q - m_p)/(q-1)} \phi^s(T_{\bj_{p}})^{-1}
\prod_{l \in L(\J_p(\bj_p))}^{[ m_p-1]}
 \Big(\sum_{|\bu| =l, \bu \succeq \bj_p} \phi^s(T_\bu)^{1-q} \mu(C_{\bu})^q 
\Big)^{1/(q-1)}\bigg]^{(q-1)/n} d\mu(\bj)  \nonumber\\
&
= 
\mu(C_{\bj_p})^{(n - m_p)/n} \big(\phi^s(T_{\bj_{p}})^{1-q}\mu(C_{\bj_p})^{q}\big)^{1/n}
\prod_{l \in  L(\J_p(\bj_p))}^{[ m_p-1]}
 \Big(\sum_{|\bu| =l, \bu \succeq \bj_p} \phi^s(T_\bu)^{1-q} \mu(C_{\bu})^q 
\Big)^{1/n}  \nonumber \\
& = 
\mu(C_{\bj_p}) ^{(n - n_p)/n}
\prod_{l \in L_r}^{[ n_p]}
 \Big(\sum_{|\bu| =l, \bu \succeq \bj_p } \ps(T_\bu)^{1-q} \mu(C_{\bu})^q\Big)^{1/n},  \nonumber
 \end{align}
on incorporating $(\phi^s(T_{\bj_{p}})^{1-q}\mu(C_{\bj_p})^{q})^{1/n}$ in the main product and noting that $m_p = n_p$. (Observe that this remains valid if $m_p =1$, in which case  the inner integral in (\ref{exp1}) is with respect to the single variable $\bi^{p}_{1}$ over the cylinder $C_{\bj_p}$, taking $\phi^s(\bi^{p}_{1})= 1$.) This establishes the inductive hypothesis (\ref{indhop}) when $r=p$.

Now assume that (\ref{indhop}) is valid for $r= k,\ldots, p$ for some $ 2\leq k \leq p$. Then 
\begin{align}
I \equiv &\int_{\bj \succeq \bj_{k-1}}\bigg[
\int_{\bigwedge(\bi^{k-1}_{1},\ldots,\bi^{k-1}_{m_{k-1}}) \in \J_{k-1}(\bj_{k-1})}
 \int_{\bigwedge(\bi^{k}_{1},\ldots,\bi^{k}_{m_{k}}) \in \J_{k}(\bj_{k})}\cdots
\int_{\bigwedge(\bi^{p}_{1},\ldots,\bi^{p}_{m_{p}}) \in \J_{p}(\bj_p)} \nonumber\\
& \quad \Big( \phi^s(\bi^{k-1}_{1},\ldots,\bi^{k-1}_{m_{k-1}})^{-1} \phi^s(T_{\bj_{k-1}})^{-1}
 \phi^s(\bi^{k}_{1},\ldots,\bi^{k}_{m_{k}})^{-1} \phi^s(T_{\bj_{k}})^{-1} \cdots
 \phi^s(\bi^{p}_{1},\ldots,\bi^{p}_{m_{p}})^{-1} \phi^s(T_{\bj_{p}})^{-1}\Big)\nonumber\\
 &\qquad \qquad \qquad d\mu(\bi^{k-1}_{1})\ldots d\mu(\bi^{k-1}_{m_{k-1}})
 d\mu(\bi^{k}_{1})\ldots d\mu(\bi^{k}_{m_k})\cdots 
 d\mu(\bi^{p}_{1})\ldots d\mu(\bi^{p}_{m_p})\bigg]^{(q-1)/n} d\mu(\bj)  \nonumber \\
  &
= \phi^s(T_{\bj_{k-1}})^{(1-q)/n}\bigg[
\int_{\bigwedge(\bi^{k-1}_{1},\ldots,\bi^{k-1}_{m_{k-1}}) \in \J_{k-1}(\bj_{k-1})}
\phi^s(\bi^{k-1}_{1},\ldots,\bi^{k-1}_{m_{k-1}})^{-1}
d\mu(\bi^{k-1}_{1})\ldots d\mu(\bi^{k-1}_{m_{k-1}})\bigg]^{(q-1)/n}\nonumber\\
&\qquad \times  \sum_ {\bj_k \succeq \bj_{k-1}}\int_{\bj \succeq \bj_{k}}\bigg[
 \int_{\bigwedge(\bi^{k}_{1},\ldots,\bi^{k}_{m_{k}}) \in \J_{k}(\bj_{k})}\cdots
\int_{\bigwedge(\bi^{p}_{1},\ldots,\bi^{p}_{m_{p}}) \in \J_{p}(\bj_p)} \Big(
 \phi^s(\bi^{k}_{1},\ldots,\bi^{k}_{m_{k}})^{-1} \phi^s(T_{\bj_{k}})^{-1} \times \cdots\nonumber\\
&\qquad \qquad \qquad\times  \phi^s(\bi^{p}_{1},\ldots,\bi^{p}_{m_{p}})^{-1} \phi^s(T_{\bj_{p}})^{-1}\Big)
 d\mu(\bi^{k}_{1})\ldots d\mu(\bi^{k}_{m_k})\cdots 
 d\mu(\bi^{p}_{1})\ldots d\mu(\bi^{p}_{m_p})\bigg]^{(q-1)/n} d\mu(\bj)  \nonumber\\ 
  &
  \leq\phi^s(T_{\bj_{k-1}})^{(1-q)/n}
 \bigg[\mu(C_{\bj_{k-1}})^{(q - m_{k-1})/n}
\prod_{l \in L({\cal J}_{k-1} (\bj_{k-1}))}^{[ m_{k-1}-1]}
 \Big(\sum_{|\bu| =l, \bu \succeq \bj_{k-1}} \phi^s(T_\bu)^{1-q} \mu(C_{\bu})^q 
\Big)^{1/n}\bigg]\nonumber\\
&\qquad \qquad \times \sum_ {\bj_k \succeq \bj_{k-1}}\mu(C_{\bj_k}) ^{(n - n_k)/n}
\prod_{l \in L_k}^{[ n_k]}
 \Big(\sum_{|\bu| =l, \bu \succeq \bj_k } \ps(T_\bu)^{1-q} \mu(C_{\bu})^q\Big)^{1/n},\label{jaykay}
 \end{align}
 where we have used Proposition \ref{integerest} to estimate the first part and the inductive hypothesis (\ref{indhop}) for the second part.
Using H\"{o}lder's inequality for each $\bj_{k-1}$:
\begin{align*}
\sum_ {\bj_k \succeq \bj_{k-1}}&\mu(C_{\bj_k}) ^{(n - n_k )/n}
\prod_{l \in L_k}^{[ n_k ]}
 \Big(\sum_{|\bu| =l, \bu \succeq \bj_k } \ps(T_\bu)^{1-q} \mu(C_{\bu})^q\Big)^{1/n} \\
& \leq\Big(\sum_ {\bj_k \succeq \bj_{k-1}}\mu(C_{\bj_k})\Big)^{(n - n_k )/n}
\prod_{l \in L_k}^{[ n_k ]}
 \Big(\sum_ {\bj_k \succeq \bj_{k-1}}\sum_{|\bu| =l, \bu \succeq \bj_k } \ps(T_\bu)^{1-q} \mu(C_{\bu})^q\Big)^{1/n} \\
& = \mu(C_{\bj_{k-1}})^{(n - n_k )/n}
\prod_{l \in L_k}^{[ n_k ]}
 \Big(\sum_{|\bu| =l, \bu \succeq \bj_{k-1} } \ps(T_\bu)^{1-q} \mu(C_{\bu})^q\Big)^{1/n}. 
\end{align*} 
Thus from (\ref{jaykay})
 \begin{align*}
I \leq \mu &(C_{\bj_{k-1}})^{(n - n_k - m_{k-1} )/n} \big(\phi^s(T_{\bj_{k-1}})^{(1-q)} \mu(C_{\bj_{k-1}})^q\big)^{1/n}
 \\
&\times \prod_{l \in L({\cal J}_{k-1} (\bj_{k-1}))}^{[ m_{k-1}-1]}
 \Big(\sum_{|\bu| =l, \bu \succeq \bj_{k-1}} \phi^s(T_\bu)^{1-q} \mu(C_{\bu})^q 
\Big)^{1/n}
\prod_{l \in L_k}^{[ n_k]}
 \Big(\sum_{|\bu| =l, \bu \succeq \bj_{k-1} } \ps(T_\bu)^{1-q} \mu(C_{\bu})^q\Big)^{1/n}\\
 = \mu &(C_{\bj_{k-1}})^{(n - n_k - m_{k-1} )/n} 
\prod_{l \in L_{k-1}}^{[ m_{k-1}+n_k]}
 \Big(\sum_{|\bu| =l, \bu \succeq \bj_{k-1} } \ps(T_\bu)^{1-q} \mu(C_{\bu})^q\Big)^{1/n} ,
 \end{align*}
which is (\ref{indhop}) with $r=k-1$, noting that $m_{k-1}+n_k = n_{k-1}$.

Finally, taking $r=1$ in (\ref{indhop}) and noting that $n_1 = n$,
\begin{align*}
\int_{\bj \succeq \bj_1 }\bigg[
\int_{\bigwedge(\bi_{1},\ldots,\bi_{m_{1}}) \in \J_{1}(\bj_1)} \cdots 
\int_{\bigwedge(\bi_{n-m_p +1},\ldots,\bi_{n}) \in \J_{p}(\bj_p)} &
 \phi^s(\bi_1,\ldots,\bi_n,\bj)^{-1} d\mu(\bi_1)\ldots d\mu(\bi_n)\bigg]^{(q-1)/n} d\mu(\bj)  \\ 
&  
\leq
\prod_{l \in L_1}^{[ n]}
 \Big(\sum_{|\bu| =l, \bu \succeq \bj_1} \ps(T_\bu)^{1-q} \mu(C_{\bu})^q\Big)^{1/n},
\end{align*}
and summing  over all $\bj_1$ at level $k_1$ and using H\"{o}lder's inequality again, gives (\ref
{qintest}).
$\Box$

\medskip

To use (\ref{fracest}) to determine when the integral in (\ref{probest2}) converges we need to bound the number of distinct join classes that have a prescribed set of levels. Let $0 \leq l_1\leq \cdots \leq l_n$ be (not necessarily distinct) levels. Write
\begin{align}
N(l_1,\ldots, l_n) = \#\big\{(k_1,&\ldots,k_p, \J_1,\ldots,\J_p): 1 \leq p \leq n, 0 \leq k_1<\cdots <k_p, \nonumber\\
&\J_r \mbox{ is a join class with root at level } k_r, \nonumber\\
&\mbox{ such that }
L(\J_1,\ldots,\J_p,k_1,\ldots,k_p) = \{l_1,\ldots, l_n\} \big\}. \label{numbdef}
\end{align}

\begin{lem}\label{count}
Let $n \in \bbbn$ and $0<\lambda<1$. Then
$$\sum_{0 \leq l_1\leq \cdots \leq l_n}N(l_1,\ldots,l_n)  \lambda^{(l_1+\cdots +l_n)/n} <\infty.$$
\end{lem}
{\it Proof.} Let  $N_0 (l_1,\ldots, l_n)$ be the total number of join classes with root $\emptyset$ (the root of the tree $I$) and levels $l_1\leq \cdots \leq l_n$.  Every join set with levels   $0 \leq l_1\leq \cdots \leq l_n\leq l_{n+1}$ may be obtained by joining a vertex at level $ l_{n+1}$ to some vertex of a join set with levels $0 \leq l_1\leq \cdots \leq l_n$ through a path in the tree $I$, and this may be done in at most $n$ inequivalent ways to within tree automorphism. Thus
$N_0 (l_1,\ldots, l_{n+1}) \leq n N_0 (l_1,\ldots, l_n)$, so since $N_0 (l_1)=1$, we have
$$N(l_1,\ldots, l_n) \leq N_0 (l_1,\ldots, l_n) \leq (n-1)!.$$

 Thus
 \begin{align}
\sum_{0 \leq l_1\leq \cdots \leq l_n}N(l_1,\ldots,l_n)  \lambda^{(l_1+\cdots +l_n)/n}
&\leq
(n-1)!\sum_{0 \leq l_1\leq \cdots \leq l_n}  \lambda^{( l_1+\cdots +l_n)/n}\nonumber\\
&\leq
(n-1)! \sum_{k=1}^{\infty}P(k)  \lambda^{k/n},\label{series}
\end{align}
where $P(k)$ is the number of distinct ways of partitioning the integer $k$ into a sum of $n$ integers
$k = l_1 + \cdots + l_{n}$ where $0 \leq l_1\leq \cdots \leq l_n$. Since $P(k)$ is polynomially bounded (trivially $P(k) \leq (k+1)^{n-1}$), (\ref{series}) converges for $0< \lambda<1$.
$\Box$

\medskip
Using Lemma \ref{count} to count the domains of integration to which Proposition \ref{fracest}  is applied leads to the main estimate.

\begin{theo}\label{mainint}
Let $s>0$ be such that 
\begin{equation}
\limsup_{k \to \infty} \frac{\log \sum_{|\bu| =k} \ps(T_\bu)^{1-q} \mu(C_{\bu})^q}{\log k} < 1.\label{condition}
\end{equation}
Then
\begin{equation*}
I \equiv \int\bigg[
 \int \cdots \int 
 \phi^s(\bi_1,\ldots,\bi_n,\bj)^{-1} d\mu(\bi_1)\ldots d\mu(\bi_n)\bigg]^{(q-1)/n} d\mu(\bj)<\infty.  \end{equation*}
\end{theo} 
{\it Proof.}  
For each $\bj$ we decompose the integral inside the square brackets as a sum of integrals taken over all $0 \leq k_1<\cdots<k_p$  and all  $m_1,\dots,m_p\geq 1$ such that $m_1+\cdots+m_p = n$, and all join classes
$\J_1(\bj_1),\ldots,\J_p(\bj_p)$ where $\J_r(\bj_r)$ has root  $\bj_r = \bj|_{k_r}$ and spread $m_r$:
\begin{align*}
&I = \\
&\int \bigg[\hspace{-0.3cm} \sum_{{\scriptsize
\begin{array}{c}
0 \leq k_1<\cdots<k_p \\
m_1+\cdots+m_p = n   \\
\J_1,\ldots,\J_p
\end{array}}}
\hspace{-0.3cm} \int_{\bigwedge(\bi_{1},\ldots,\bi_{m_{1}}) \in \J_{1}(\bj_1)}\hspace{-0.5cm} \cdots
\int_{\bigwedge(\bi_{n-m_p +1},\ldots,\bi_{n}) \in \J_{p}(\bj_p)}
\hspace{-0.5cm} \phi^s(\bi_1,\ldots,\bi_n,\bj)^{-1} d\mu(\bi_1)\ldots d\mu(\bi_n)\bigg]^{(q-1)/n} \hspace{-0.5cm} d\mu(\bj)  \\ 
&
\leq\hspace{-0.5cm}\sum_{{\scriptsize
\begin{array}{c}
0 \leq k_1<\cdots<k_p \\
m_1+\cdots+m_p = n   \\
\J_1,\ldots,\J_p
\end{array}}}
\hspace{-0.3cm} \int \bigg[\int_{\bigwedge(\bi_{1},\ldots,\bi_{m_{1}}) \in \J_{1}(\bj_1)} \hspace{-0.5cm}\cdots
\int_{\bigwedge(\bi_{n-m_p +1},\ldots,\bi_{n}) \in \J_{p}(\bj_p)}
\hspace{-0.5cm} \phi^s(\bi_1,\ldots,\bi_n,\bj)^{-1} d\mu(\bi_1)\ldots d\mu(\bi_n)\bigg]^{(q-1)/n} d\mu(\bj)  \\ 
&
\leq\sum_{{\scriptsize
\begin{array}{c}
0 \leq k_1<\cdots<k_p \\
m_1+\cdots+m_p = n   \\
\J_1,\ldots,\J_p
\end{array}}}
 \prod_{l \in L}^{[n]}
 \Big(\sum_{\bu \in l} \ps(T_\bu)^{1-q} \mu(C_{\bu})^q\Big)^{1/ n} 
\end{align*}
where the product is over the set of levels
$L = \{L(\J_r),\ldots, L(\J_p),k_r,\ldots,k_p\}$ counted with repetitions, and we have used Minkowski's inequality and (\ref{qintest}).

Condition (\ref{condition}) implies that 
$\sum_{|\bu| =k} \ps(\bu)^{1-q} \mu(C_{\bu})^q \leq c \lambda^k$ for all $k$, for some $c>0$ and some $ \lambda <1$. Thus
\begin{align*}
I &
\leq\sum_{{\scriptsize
\begin{array}{c}
0 \leq k_1<\cdots<k_p \\
m_1+\cdots+m_p = n   \\
\J_1,\ldots,\J_p
\end{array}}}
 \prod_{l \in L}^{[n]}
(c\lambda^l)^{1/ n} \\
& 
\leq
 c\sum_{l_1,\ldots,l_n}N(l_1,\ldots,l_n)  \lambda^{(l_1+\cdots +l_n)/n} <\infty
 \end{align*}
using Lemma \ref{count}.
$\Box$

\section{Main results}
\setcounter{equation}{0}
\setcounter{theo}{0}

We now put together the estimates from the two preceding sections to obtain an almost sure lower bound for the lower $q$-dimension of measures on almost self-affine sets, which coincides with the upper bound of Corollary \ref{corub}. We then consider the special cases where the underlying measure $\mu$ is a Bernoulli measure or a Gibbs measure on $I_\infty$ when it turns out that the lower and upper $q$-dimensions coincide almost surely and there are further natural expressions for this common value. Recall the expressions for $d_{q}^{-}$ and $d_{q}^{+}$ given by (\ref{basicb1})-(\ref{seriessum}).

\begin{theo}\label{thm8.1}
Let $T_{i}$ be linear contractions on $\Esp$ with
$\| T_{i} \| < 1 $ 
for all $i = 1,\ldots,m$. 
Let $\mu$ be a finite Borel measure on $\ii$ and let $\mu^{\w}$ be the measure defined by
$(\ref{mesdef})$ in the random model described in Section \ref{rand}.  For $q>1$, we have that, for 
almost all $\w$,
$$
\lqd (\mu^\w) = \min \{ d_{q}^{-}(T_{1},\ldots,T_{m};\mu),N \},
$$
where
$$
d_{q}^{-}(T_{1},\ldots,T_{m};\mu)=\sup \{s :\sum_{k=0}^{\infty} \sum_{\bi \in I_{k}}  
 \svs {T_{\bi}} ^{1-q} 
\mu (C_{\bi})^{q} <\infty \}.
$$
\end{theo} 

\noindent{\it Proof.} 
Corollary \ref{corub} gives that $\lqd (\mu^\w) \leq \min \{ d_{q}^{-}(T_{1},\ldots,T_{m};\mu),N \}$. 

Taking  non-integral $s$ and $s_1$ such that $0<s< s_1 <d_q^-$, it follows from (\ref{basicb1}) (noting, as before, that 
$\phi^{s}(T_\bi) \geq   \alpha_{+}^{-k(s_1-s)} \phi^{s_1}(T_\bi)$ if
$|\bi|= k$) that 
$$ \sum_{|\bi| =k} \ps(T_\bi)^{1-q} \mu(C_{\bi})^q \leq c_2 \lambda^k
$$
for some $\lambda <1$, so condition (\ref{condition}) is satisfied  It follows from Theorem \ref{mainint} and Proposition
\ref{propexp} that, for all $0<s_2 <s$,
$$
\E \int r^{s_2(1-q)} \mu^{\omega}(B(x,  r))^{q-1}  d\mu^{\omega}(x) 
\leq  M r^{(s-s_2)(q-1)}
$$
for all sufficiently small $r$, for some  $M<\infty$. For any $0<\rho<1$, the Borel-Cantelli Lemma implies that almost surely the sequence
$$\int (\rho^k)^{s_2(1-q)} \mu^{\omega}(B(x,  \rho^k)^{q-1}  d\mu^{\omega}(x) \quad (k=1,2,\ldots)$$
converges to $0$, so, since the asymptotic behaviour of the multifractal integrals is controlled by their values on any such sequence of $r = \rho^k$, we conclude that 
$$\lim_{r \to 0} \int r^{s_2(1-q)} \mu^{\omega}(B(x,  r))^{q-1}  d\mu^{\omega}(x) =0$$
almost surely. This is true for all $s_2< d_q^-$, so $\lqd (\mu^\w) \geq \min \{ d_{q}^{-}(T_{1},\ldots,T_{m};\mu),N \}$ almost surely.
$\Box$

\medskip

We now specialise to Bernoulli measures on almost self-affine sets, which might be termed `almost self-affine measures'.
Let $p_{1},\ldots ,p_{m}$ be `probabilities', 
with $p_{i} > 0$ for all $i$ and  $\sum_{i=1}^{m}p_{i} = 1$.
We may define a self-similar Borel 
measure $\mu$ on $\ii$ by setting 
\begin{eqnarray}
\mu (C_{\bi}) = p_{\bi} \equiv  p_{i_{1}} p_{i_{2}}\ldots 
p_{i_{k}},\label{probs}
\end{eqnarray}
on the cylinders $C_{\bi}$, where  ${\bf i} = (i_{1}, \ldots , i_{k})$, 
and extending to general 
Borel and measurable subsets of $\ii$ in the usual way. (The measure $\mu$ may be thought of as an invariant measure on the code space $\ii$ under the shift map.) We refer to the measures $\mu^\w$ on the almost self-affine sets $E^{\w}$ as {\it Bernoulli measures} on $E^{\w}$ or {\it almost self-affine measures}.

\begin{lem}\label{lem8.2}
Let $\mu$ be defined by $(\ref{probs})$.  
For all $q > 1$,
the limit
\begin{eqnarray}
\lim_{k \rightarrow \infty}\Big(\sum_{\bi \in I_{k}}
\svf {s} {\ti} ^{1-q} \mu(\ci)^{q}\Big) ^{1/k}  
=\lim_{k \rightarrow \infty}\Big(\sum_{\bi \in I_{k}}
\svf {s} {\ti} ^{1-q}p_{\bi} ^{q}\Big) ^{1/k}  
\label{3.9}
\end{eqnarray}
exists  for all $s>0$ and is strictly increasing in $s$. In particular there is a unique number $d_{q} \equiv d_{q}(T_{1},\ldots,T_{m};\mu)$ such that
\begin{eqnarray}
\lim_{k \rightarrow \infty}\Big(\sum_{\bi \in I_{k}}
\svf {d_{q}} {\ti} ^{1-q} \mu(\ci)^{q}\Big) ^{1/k} = 1, 
\label{3.9a}
\end{eqnarray}
and, moreover, $d_{q}= d_{q}^{-}=d_{q}^{+}$.
\end{lem}

\noindent{\it Proof. }
With $\mu$ a Bernoulli measure, it follows from (\ref{probs}) and  (\ref{2.2}) that $\big\{ \sum_{\bi \in I_{k}}
\svf {s} {\ti} ^{1-q} \mu(\ci)^{q}\big\}_{k=0}^{\infty}$ is a supermultiplicative sequence, so by the standard property of such sequences, the limit (\ref{3.9}) exists. Monotonicity, and thus the existence of a unique $d_{q}$ satisfying (\ref{3.9a}), follows from (\ref{phicont}). That  $d_{q}^- = d_{q}$ follows from (\ref{lgqdim}). The argument of \cite[Proposition 6.1]{Fa5} establishes that $d_{q}^+ = d_{q}$.
$\Box$
\medskip

Our result for almost self-affine measures now follows easily.

\begin{cor}\label{cor8.3}
{\em (Almost self-affine measures)} Let $T_{i}$ be linear contractions on $\Esp$ with
$\| T_{i} \| < 1 $ 
for all $i = 1,\ldots,m$. 
Let $\mu$ be the Bernoulli measure on $\ii$ given by $(\ref{probs})$ and let $\mu^{\w}$ be the almost self-affine measure defined by
$(\ref{mesdef})$ on  the random set $E^{\w}$. Then, for 
almost all $\w$, 
$$
\qd(\mu^\w)=  \lqd (\mu^\w) = \uqd (\mu^\w) = \min \{ d_{q}(T_{1},\ldots,T_{m};\mu),N \},
$$ 
for all $q>1$, where $d_{q}(T_{1},\ldots,T_{m};\mu)$ is the unique positive number satisfying
$$\lim_{k \rightarrow \infty}\Big(\sum_{\bi \in I_{k}}
\svf {d_{q}} {\ti} ^{1-q} \mu(\ci)^{q}\Big) ^{1/k} = 1.
$$
\end{cor} 

\noindent{\it Proof.} 
By Theorem \ref{thm8.1}  and Corollary \ref{corub} and we have that, almost surely,  
$$ \min\{d_{q}^-,N\}  =   \lqd (\mu^\w)\leq \uqd (\mu^\w) \leq \min\{d_{q}^+,N\} ,$$ 
so the conclusion follows from Lemma \ref{lem8.2}.
$\Box$.

\medskip

As might be anticipated, if $\mu$ is a Gibbs measure on $\ii$ we get similar results to those 
for $\mu$ a Bernoulli measure.  Let $\sigma$ be the shift map on $\ii$, 
so $\sigma (i_{1},i_{2},\ldots) = (i_{2},i_{3},\ldots)$.  For 
$f: \ii \rightarrow \bbbr$ we define the sums 
\begin{equation}
	S_{k}f(\bi) = \sum_{j=0}^{k-1} f(\sigma^{j}(\bi)), \label{sum}
\end{equation}
where $\bi = (i_{1},i_{2},\ldots) \in I_\infty$, and $\sigma^{j}$ is the $j$th 
iterate of $\sigma$.  A Borel probability measure 
$\mu$ on $\ii$ is a {\it Gibbs measure} on $I_{\infty}$ if there 
exists a continuous $f: \ii \rightarrow \bbbr$, a number $P(f)$ 
termed the {\it pressure} of $f$, and $a>0$, such that for all $k$ and all 
$\bi = (i_{1},\ldots,i_{k}) \in I_{k}$ we have 
\begin{equation}
a \leq \frac{\mu(C_{\bi})}{\exp(-kP(f) + S_{k}f(\bi))} \leq a^{-1}. \label{gibbs}
\end{equation}
Thus the pressure is given by
\begin{equation*}
P(f) = \lim_{k \rightarrow \infty}\frac{1}{k}\log
\sum_{\bi \in I_{k}} \exp (S_{k}f(\bi)). 
\end{equation*}
By a standard result from the thermodynamic formalism, see for example \cite{Fa3,P}, if $f$ satisfies an 
$\epsilon$-H\"{o}lder condition of the form $|f(\bi)-f(\bj)| \leq c 
d(\bi,\bj)^{\epsilon}$ for all $\bi,\bj \in \ii$ for some $\epsilon >0$, 
then there exists an invariant Gibbs measure $\mu$ 
satisfying (\ref{gibbs}) for some $P(f)$.

From (\ref{sum})
\begin{equation*}
	S_{k+l}f(\bi) = S_{k}f(\bi) + S_{l}f(\sigma^{k}\bi) 
\end{equation*}	
for $k,l = 1,2,\ldots$, so from (\ref{gibbs})
\begin{equation}
a^{3} \leq \frac{\mu(C_{\bi,\bj})}{\mu(C_{\bi})\mu(C_{\bj})} \leq 
a^{-3}. \label{mesprod}
\end{equation}
for all $\bi,\bj \in I$.  This inequality leads to analogues of Lemma \ref{lem8.2} and Corollary \ref{cor8.3} for Gibbs measures.

\begin{lem}\label{lem8.4}
The conclusions of Lemma $\ref{lem8.2}$ hold if $\mu$ is a Gibbs measure satisfying $(\ref{gibbs})$.
\end{lem}

\noindent{\it Proof. }
It follows from (\ref{mesprod}) and (\ref{2.2}) that  $\big\{ a^{3q} \sum_{\bi \in I_{k}}
\svf {s} {\ti} ^{1-q} \mu(\ci)^{q}\big\}_{k=0}^{\infty}$ is a supermultiplicative sequence, so again the limits (\ref{3.9}) exist. The other conclusions follow just as in Lemma \ref{lem8.2}, see also \cite[Proposition 7.1]{Fa5}.
$\Box$

\medskip

The result for Gibbs measures now follows.

\begin{cor}\label{cor8.5}
{\em (Gibbs measures)} Let $T_{i}$ be linear contractions on $\Esp$ with
$\| T_{i} \| < 1 $ 
for all $i = 1,\ldots,m$. 
Let $\mu$ be a Gibbs measure on $\ii$ satisfying $(\ref{gibbs})$ and let $\mu^{\w}$ be the almost self-affine measure on  $E^{\w}$ defined by
$(\ref{mesdef})$.  For  the random model described in Section \ref{rand}, for 
almost all $\w$ we have
$$
\qd(\mu^\w)=  \lqd (\mu^\w) = \uqd (\mu^\w) = \min \{ d_{q}(T_{1},\ldots,T_{m};\mu),N \},
$$ 
for all $q>1$, where $d_{q}(T_{1},\ldots,T_{m};\mu)$ is the unique positive number satisfying
$$\lim_{k \rightarrow \infty}\Big(\sum_{\bi \in I_{k}}
\svf {d_{q}} {\ti} ^{1-q} \mu(\ci)^{q}\Big) ^{1/k} = 1.
$$
\end{cor} 

\noindent{\it Proof.} 
This is precisely as in   Corollary \ref{cor8.3}, using Lemma \ref{lem8.4} rather than Lemma \ref{lem8.2}.
$\Box$.

\medskip

The expression 
$\lim_{k \rightarrow \infty}\left(\sum_{\bi \in I_{k}}
\svf {s} {\ti} ^{1-q} \mu(\ci)^{q}\right) ^{1/k} = 1$, 
that gives the generalised dimensions for Gibbs measures satisfying $(\ref{gibbs})$, may be regarded as (the exponential of a)  pressure expression in the context of a subadditive or generalised thermodynamic formalism, see \cite{Bar, Fa4}.  With an appropriate definition of generalised pressure $P(\{g_k\})$ for a subadditive family of functions $\{g_k\}$, the number $d_{q}$ is the unique value of $s$ such that 
$$
P\big(\big\{(1-q) \log \phi^{s}(T_{\cdot |k})+ q(S_{k}f(\cdot ) - k P(f))\big\}_{k}\big) = 0 
$$
or equivalently
$$
P\big(\big\{(1-q) \log \phi^{s}(T_{\cdot |k})+ q\mu(C_{\cdot |k}\big\}_{k}\big) = 0,
$$
see \cite{Fa5} for more details. 

\section{Further remarks}
\setcounter{equation}{0}
\setcounter{theo}{0}

(1) The numbers $d_{q}$ can have discontinuous derivatives 
at values of $q$ for which $d_{q}$ is an integer,  since for $s$ non-integral
$$\frac{d\svs T }{ds}= \svs T \log \alpha_{j},$$
where $j$ is the integer such that $j-1 < s < j$, which
is discontiouous at the integer $s$ if $\alpha_{s}(T) > 
\alpha_{s+1}(T)$.  
Thus the  $q$-dimensions of measures on almost self-affine sets 
typically exhibit phase transitions, that is have discontinuous derivatives with respect to $q$.  

For a simple example, let $T:\Esp \rightarrow \Esp$ be a 
self-adjoint linear mapping with distinct singular values and with $| T| < 1$, 
let  $T_{1}= \cdots =  T_{m} = T$, and let $\mu$ 
be the Bernoulli measures defined by (\ref{probs}).  It is easily checked that $d_{q}$ is defined by the requirement that
$\svf {d_{q}} {T}  
= \left(\sum_{i=1}^{m}p_{i}^{q}\right)^{1/(q-1)}$, so almost surely, the generalised dimensions 
will not be differentiable at values of $q$ where  $\qd(\mu^\w)$ takes integer values.

(2) The conditions on the distribution of random displacements $\w_\bi$ stated in Section 5 can be weakened considerably with the main results of Section 8 still holding. The arguments go through unchanged if the random vectors $\w_\bi \in D$ are independent with uniformly bounded density  -- identical distribution is not essential.

(3)  It is natural to ask under what other conditions the `generic'  formula (\ref{3.9a}) gives the generalised dimensions. In particular, what can be said for measures on self-affine sets rather than almost self-affine sets? Whilst the formula holds for almost all strictly self-affine sets (with respect to translates  $a_1,\ldots,a_m \in \Esp$) if $1<q\leq 2$,  more randomness seems to be unavoidable if $q>2$.

Finding generic expressions for $q$-dimensions of self-affine-like sets if $0<q<1$  needs a different approach.   Whilst  $d_q$ (defined with an infimum in (\ref{affdim})) provides an upper bound it seems awkward to show that this is the generic value.  In view of examples related to the projection of measures 
where the `natural' formulae for $q$-dimensions fail for  $0<q<1$, see \cite{HK},  a generic formula might well be more subtle.

(4) It would be of interest to develop a   `fine' multifractal analysis  of measures on
(almost) self-affine sets, and to find generic forms of the
(Hausdorff)
multifractal spectrum $f_{\rm H} (\alpha)$ of $\mu^\w$, 
that is 
$$f_{\rm H} (\alpha) = \dim_{H} \{ x \in \Esp: \lim_{r \rightarrow 0} 
\log \mu^\w (B(x,r)) /\log r = \alpha \},$$
where $\dim_{H}$ denotes Hausdorff dimension.
From general results on coarse and fine 
multifractal theory and their relationships, see, for example, \cite[Chapter 11]{Fa3}.


\begin{thebibliography}{abc-12}
 \bibitem{BM} J. Barral and M. Mensi. Gibbs measures on self-affine Sierpinski carpets and their singularity
spectrum. {\em Ergod. Th. Dynam. Sys. }{\bf 27} (2007) 1419-1443.
    \bibitem{Bar} L.M. Barreira. A non-additive thermodynamic 
    formalism and applications to dimension theory of hyperbolic 
    dynamical systems, {\em Ergod. Th. Dynam. Sys. }{\bf 16} 
    (1996)~871-927.
 \bibitem{Bed} T. Bedford. {\em Crinkly curves, Markov partitions and box 
dimensions in self-similar sets}\/ (PhD thesis, University of 
Warwick, 1984).
 \bibitem{Edg}G.A. Edgar. Fractal dimension of self-affine sets: some examples. {\em Rend. Circ. Mat.
Palermo (2) Suppl. }{\bf 28}(1988)  341-358.
   \bibitem{Fa1} K.J.~Falconer. The Hausdorff dimension of 
   self-affine fractals, {\em Math. Proc. Cambridge Philos. 
Soc. }{\bf 103} (1988)~339-350.
   \bibitem{Fa2} K.J.~Falconer. The dimension of 
   self-affine fractals II, {\em Math. Proc. Cambridge Philos. 
Soc. }{\bf 111} (1992)~169-179.
\bibitem{Fa4} K.J.~Falconer. Bounded distortion and dimension for 
 non-conformal repellers, {\em Math. Proc. Cambridge Philos. 
Soc. }{\bf 115} (1994)~315-334.
    \bibitem{Fa3} K.J.~Falconer. {\em Techniques in Fractal Geometry}
    \/ (John Wiley, 1997).
  \bibitem{Fa5} K.J.~Falconer. Generalized dimensions of measures on
   self-affine sets, {\em Nonlinearity }{\bf 12} (1999)~877-891.
  \bibitem{Fa} K.J.~Falconer. {\em Fractal Geometry---Mathematical 
Foundations and Applications}, 2nd Ed. (John Wiley, 2003).

\bibitem{Gr} P. Grassberger. Generalised dimension of strange 
   attractors, {\em Phys. Rev. Lett. A }{\bf 97} (1983)~227-230.
 \bibitem{Har} D. Harte. {\em Multifractals: Theory and Applications}
    \/ (Chapman and Hall, 2001).
  \bibitem{HL} I. Heuter and S. Lalley. Falconer's formula for the 
    Hausdorff dimension of a self-affine set in $\bbbr^{2}$,  
    {\em Ergod. Th. Dynam. Sys. }{\bf 15} 
    (1995)~77-97.
  \bibitem{HK} V.Y. Hunt and B.R. Kaloshin.  How projections affect 
   the dimensions of fractal measures, {\em Nonlinearity }{\bf 10} 
   (1997)~1031-1046.
  \bibitem{JPS}  T. Jordan, M. Pollicott and K. Simon. Hausdorff dimension for randomly
   perturbed self affine attractors, {\em Commun. Math. Phys. }{\bf 270} 
    (2007)~519-544.
\bibitem{K} J. King. The singularity spectrum for general 
Sierpinski carpets, {\em Adv. Math. }{\bf 116} (1995)~1-8.

\bibitem{KS} A. K\"{a}enm\"{a}ki and P. Shmirkin. Overlapping self-affine sets of Kakeya type, 
{\em Ergod. Th. Dynam. Sys. }{\bf 29} (2009) 941-965.

\bibitem{Lau} K.-S. Lau. Self-similarity, $L^{p}$-spectrum and 
   multifractal formalism, in  Fractal Geometry and Stochastics,
   Eds. C. Bandt, S. Graf and M. Z\"{a}hle, {\em Progress in Probability } {\bf 37}  
   55-90 (Birkh\"{a}user, 
   1995).
\bibitem{Man} B. Mandelbrot. Negative fractal dimensions and 
multifractals, {\em Physica A }{\bf 163} (1990)~306-315.
   \bibitem{McM} C. McMullen. The Hausdorff dimension of 
   general Sierpi\'{n}ski carpets, {\em Nagoya Math. 
   J. }{\bf 96} (1984)~1-9.
  \bibitem{O} L. Olsen. Self-affine multifractal Sierpinski sponges in 
  $\bbbr^{d}$, {\em Pacific J. Math. }{\bf 183} (1998)~143-199.  
 \bibitem{PSol} 
Y. Peres and B. Solomyak, Problems on self-similar sets and self-affine sets: an update, in Fractal Geometry and Stochastics II,  Eds. C. Bandt, S. Graf, and M. ZŠhle, 
{\em Progress in Probability } {\bf 46} 95-106, Birkh\"{a}user, 2000.
\bibitem{P} Y. B. Pesin. {\em Dimension theory in dynamical 
systems} (University of Chicago Press, 1997). 
	\bibitem{Sol} B. Solomyak. Measure and dimensions for some 
  fractal families, {\em  Math. Proc. Cambridge Philos. Soc. } {\bf 124} (1998) 531Ð546.


\end{thebibliography}
\end{document}